\documentclass[oneside]{amsart} 

\usepackage[percent]{overpic}

\usepackage{paralist}
\usepackage{graphics} 
\usepackage{epsfig} 
\usepackage{xspace}
\usepackage{marginnote} 

\usepackage{array}

\usepackage{appendix}
\usepackage[english]{babel} 
\usepackage[T1]{fontenc} 
\usepackage[utf8]{inputenc} 
\usepackage{geometry}
\usepackage{graphicx}
\usepackage{tabto}
\usepackage{amsmath}
\usepackage{calligra}
\usepackage{amsthm}
\usepackage{bbm}
\usepackage{cancel}
\usepackage{mathrsfs}
\usepackage{mathtools}
\usepackage{version}
\usepackage{arydshln}
\usepackage[math]{cellspace}

\usepackage{multicol}

\usepackage{bm}
\usepackage{subfig}

\usepackage{tikz}

\usepackage{amssymb}
\usepackage{stmaryrd}
\newcommand{\numberset}{\mathbb}

\newcommand{\R}{\numberset{R}}
\newcommand{\F}{\numberset{F}}
\newcommand{\I}{\numberset{I}}
\newcommand{\PP}{\numberset{P}}

\DeclareFontFamily{U}{matha}{\hyphenchar\font45}
\DeclareFontShape{U}{matha}{m}{n}{
	<-6> matha5 <6-7> matha6 <7-8> matha7
	<8-9> matha8 <9-10> matha9
	<10-12> matha10 <12-> matha12
}{}
\DeclareSymbolFont{matha}{U}{matha}{m}{n}
\DeclareFontFamily{U}{mathx}{\hyphenchar\font45}
\DeclareFontShape{U}{mathx}{m}{n}{
	<-6> mathx5 <6-7> mathx6 <7-8> mathx7
	<8-9> mathx8 <9-10> mathx9
	<10-12> mathx10 <12-> mathx12
}{}
\DeclareSymbolFont{mathx}{U}{mathx}{m}{n}
\DeclareMathDelimiter{\vvvert} {0}{matha}{"7E}{mathx}{"17}%
\newcommand{\norm}[1]{\left\lVert#1\right\rVert}

\newcommand{\jump}[1]{\left\llbracket #1 \right\rrbracket}

\newtheorem{teo}{Theorem}[section]

\newtheorem{pro}[teo]{Problem}

\newtheorem{prop}[teo]{Proposition}
\newtheorem{rem}[teo]{Remark}

\let\div\undefined\DeclareMathOperator{\div}{div} 
\let\curl\undefined\DeclareMathOperator{\curl}{curl} 
\DeclareMathOperator*{\Grad}{\boldsymbol\nabla}
\DeclareMathOperator*{\Grads}{{\underline{\boldsymbol{\varepsilon}}}}
\DeclareMathOperator{\grad}{\nabla} 
\DeclareMathOperator{\grads}{\boldsymbol\nabla_s}
\newcommand{\derivative}[2]{\frac{\partial #1}{\partial #2}}
\newcommand\Huo{\mathbf{H}^1_0(\Omega)}
\newcommand\Ldo{L^2_0(\Omega)} 
\newcommand\Hub{\mathbf{H}^1(\B)} 
\newcommand\Hubd{(\Hub)^\prime} 
\newcommand\LdB{L^2(\B)} 
\newcommand\LdO{L^2(\Omega)} 
\newcommand\LdBd{\mathbf{L}^2(\B)} 

\newcommand\HuOd{H^1(\Omega_2)}

\newcommand{\LdOdd}{L^2(\Omega_2)}
\newcommand\Oft{\Omega^f_t} 
\newcommand\Ost{\Omega^s_t} 
\newcommand\Of{\Omega^f} 
\newcommand\Os{\Omega^s} 
\newcommand\B{\mathcal B} 
\renewcommand\u{\mathbf{u}} 
\renewcommand\v{\mathbf{v}} 

\renewcommand\c{\mathbf{c}} 
\newcommand\n{\mathbf{n}} 
\newcommand\w{\mathbf{w}} 
\newcommand\x{\mathbf{x}} 
\newcommand\s{\mathbf{s}} 
\renewcommand\S{\mathbf{S}} 

\newcommand\X{\mathbf{X}} 
\newcommand\Y{\mathbf{Y}} 

\newcommand\LL{{\boldsymbol{\Lambda}}} 
\newcommand\Vline{\mathbf{V}} 
\newcommand\Sline{\S} 

\newcommand\ssigma{\boldsymbol\sigma} 
\newcommand\llambda{\boldsymbol\lambda} 
\newcommand\mmu{\boldsymbol\mu} 
\newcommand\dt{\Delta t} 
\newcommand\dr{\delta\rho} 
\newcommand\T{\mathcal{T}} 
\newcommand\quadnode{\mathbf{p}}
\newcommand\quadweigth{\omega}

\newcommand\Pcal{\mathcal{P}}

\newcommand\Bbub{\mathfrak{B}}

\newcommand{\Amatr}{\mathsf{A}}
\newcommand{\Bmatr}{\mathsf{B}}
\newcommand{\Cmatr}{\mathsf{C}}
\newcommand{\Kmatr}{\mathsf{K}}
\newcommand{\Mmatr}{\mathsf{M}}
\newcommand{\Zmatr}{\mathsf{0}}

\newcommand\Ed{K_2}
\newcommand\E{K}

\newcommand\vrd{v_{|\Omega_2}}

\newcommand\urd{u_{|\Omega_2}}
\newcommand\vhrd{v_{h|\Omega_2}}
\newcommand\uhrd{u_{h|\Omega_2}}
\newcommand\uth{u_{2,h}}
\newcommand\vth{v_{2,h}}
\newcommand\Vth{V_{2,h}}
\newcommand\Vt{H^1(\Omega_2)}
\newcommand\Vtd{(\Vt)^\prime}
\newcommand\Huzo{H^1_0(\Omega)}

\newcommand\Poly{\mathcal{P}}
\newcommand{\Qoly}{\mathcal{Q}}
\newcommand{\Boly}{\mathfrak{B}}

\newcommand\edge{e}
\newcommand\edged{\edge_2}

\newcommand{\dd}{\mathsf{dd}}
\newcommand{\md}{\mathsf{md}}
\newcommand{\mm}{\mathsf{mm}}

\usepackage{listings}
\renewcommand\lg{\begin{color}{black}}
	\newcommand\gl{\end{color}}
\newcommand\blue{\begin{color}{black}}
	\newcommand\noblue{\end{color}}

\newcommand\red{\begin{color}{black}}
	\newcommand\nored{\end{color}}

\newcommand\DB{\begin{color}{black}}
	\newcommand\BD{\end{color}}

\newcommand\fc{\begin{color}{black}}
	\newcommand\cf{\end{color}}

\newcommand\fcc{\begin{color}{black}}
	\newcommand\cff{\end{color}}

\definecolor{UBCblue}{rgb}{0.2,0.4,0.7}
\definecolor{beaver}{rgb}{0.62, 0.51, 0.44}
\definecolor{airforceblue}{rgb}{0.36, 0.54, 0.66}
\definecolor{cadmiumgreen}{rgb}{0.0, 0.42, 0.24}
\definecolor{bazaar}{rgb}{0.6, 0.47, 0.48}
\definecolor{tealblue}{rgb}{0.21, 0.46, 0.53}
\definecolor{amber}{rgb}{1.0, 0.75, 0.0}
\definecolor{bronze}{rgb}{0.8, 0.5, 0.2}
\definecolor{cadetblue}{rgb}{0.37, 0.62, 0.63}
\definecolor{coralred}{rgb}{1.0, 0.25, 0.25}
\definecolor{burlywood}{rgb}{0.87, 0.72, 0.53}

\begin{document}
	
	\title[]{Advances on finite element discretization\\of fluid-structure interaction problems}
	
	\author{Najwa Alshehri}
	\address{Department of General Studies, Jubail Industrial College, Jubail Industrial City 25718, Saudi Arabia}
	\email{shehrin@rcjy.edu.sa}
	
	\author[Daniele Boffi]{Daniele Boffi}
	\address{Computer, electrical and mathematical sciences and engineering division, King Abdullah University of Science and Technology, Thuwal 23955, Saudi Arabia and Dipartimento di Matematica \textquoteleft F. Casorati\textquoteright, Universit\`a degli Studi di Pavia, via Ferrata 5, 27100, Pavia, Italy}
	\email{daniele.boffi@kaust.edu.sa}
	\urladdr{kaust.edu.sa/en/study/faculty/daniele-boffi}
	
	\author{Fabio Credali}
	\address{Computer, electrical and mathematical sciences and engineering division, King Abdullah University of Science and Technology, Thuwal 23955, Saudi Arabia}
	\email{fabio.credali@kaust.edu.sa}
	
	\author{Lucia Gastaldi}
	\address{Dipartimento di Ingegneria Civile, Architettura, Territorio, Ambiente e di Matematica, Universit\`a degli Studi di Brescia, via Branze 43, 25123, Brescia, Italy}
	\email{lucia.gastaldi@unibs.it}
	\urladdr{lucia-gastaldi.unibs.it}
	
	\maketitle
	
	\begin{abstract}
		 We review the main features of an unfitted finite element method for interface and fluid-structure interaction problems based on a distributed Lagrange multiplier in the spirit of the fictitious domain approach. We recall our theoretical findings concerning well-posedness, stability, and convergence of the numerical schemes, and discuss the related computational challenges. In the case of elliptic interface problems, we also present a posteriori error estimates.
	\end{abstract}

	\section*{Introduction}
	
	The mathematical description of many physical and engineering problems requires partial differential equations involving different operators or discontinuous coefficients so that the computational domain is partitioned into several (possibly time dependent) regions. The coupling between different models is enforced at the interface via suitable transmission conditions. Among the applications, we mention, \fc for instance, thermo-mechanical problems~\cite{thermo-mech-1,thermo-mech-2}, cardiac simulations~\cite{peskin1972flow,balzani2016}, computational geo-science~\cite{geo-2,geo-1}\cf, and many others.
	
	When designing effective numerical methods for the mentioned problems, an accurate representation of the underlying geometry is needed. \textit{Fitted approaches} are based on a mesh which conforms to the interface and provide accurate solutions~\cite{barrett-interface,hou-interface,huynh-interface,mu-interface,chen-interface}. However, when applied to evolution problems, computational issues may arise. For instance, the Arbitrary Lagrangian--Eulerian approach~\cite{hirt1974arbitrary,donea1982arbitrary,improved-ale}, appealing for fluid-structure interaction problems, needs an update of the mesh at each time step. In case of large displacements or deformations, the update may break shape regularity, thus requiring the generation of a new mesh, a costly operation in general. Moreover, when fluid and solid have similar or equal densities, the added mass effect may affect the performance of the method~\cite{causin}. Such phenomenon may lead to the unconditional instability of the scheme regardless of the physical parameters. The added mass effect can be reduced by an appropriate treatment of transmission conditions, as described in~\cite{deparis,badia}. 
	
	Due to the critical aspects of fitted approaches, several \textit{unfitted methods} have been introduced during the past decades. In this case, meshes are generated without taking care of the interface, which is allowed to cross elements. We mention, for instance, the Immersed Boundary Method~\cite{peskin1972flow,peskin2002immersed}, where the interface position is tracked by Dirac delta functions, and the level set method~\cite{SUSSMAN1994146,chang1996level}, where the interface corresponds to the zero level set of a certain function. Other techniques such as Nitsche-XFEM~\cite{alauzet2016nitsche}, Cut-FEM~\cite{burman2015cutfem,cutStokes,HANSBO201490}, Finite Cell methods~\cite{finitecell1,finitecell2,finitecell3}, and Ghost-FEM~\cite{ghost} are based on the use of penalty terms. The fictitious domain approach has been introduced to model the particulate flow as described in~\cite{glowinski1997lagrange,glowinski2001fictitious}, while the application to more general cases can be found in~\cite{yu2005dlm,fatboundary}. \fc Another paradigm for tackling moving boundaries is given by \textit{diffuse interface approaches}~\cite{diffuse-int-1,diffuse-int-2,diffuse-int-3}, where quantities localized on the interface are distributed throughout a thick interfacial region.\cf
	
	In~\cite{auricchio2015fictitious,2015}, we introduced a new fictitious domain formulation with distributed Lagrange multiplier (FD-DLM) for the discretization of elliptic interface and fluid-structure interaction problems. Indeed, these problems share the common feature that the domain is divided into two regions separated by an interface. In the case of elliptic interface problems this is due to the presence of discontinuous coefficients, while fluid-structure interactions involve different types of equations, corresponding to constitutive models for fluid and solid. \fcc This survey paper aims to collect, \cff\fc in a unified framework\cf\fcc, the main properties of this approach, including the well-posedness and stability of the solution, convergence of the discrete scheme and computational aspects. \cff
	
	The main idea of our approach is to extend one domain into the other. The two domains are then discretized by fixed independent meshes. In particular, in case of fluid-structure interactions, the fluid dynamics is described in Eulerian framework on the extended domain, while the deformation of the solid body is represented in Lagrangian setting by considering a reference domain which is mapped, at each time instant, into the actual configuration. In order to impose that the solution in the fictitious region coincides with the solution in the immersed domain, we consider a Lagrange multiplier and a coupling term is added to the model. We emphasize that all computations are done on fixed domains, so that meshes are generated only once.
	
	The paper is organized into two parts. In Section~\ref{sec:interface}, we focus on elliptic interface problems and discuss in some details the main features of our approach, which will be adapted later on for fluid-structure interactions. After presenting the continuous and discrete problems, we describe how to deal with the coupling term~\cite{boffi2022interface,BCG24}. Indeed, the construction of its finite element counterpart requires integration over non-matching grids, which is a challenging aspect of several immersed methods~\cite{maday2002influence,farrell2011conservative,massing2013efficient,bvrezina2017fast,fromm2023interpolation}. \fc Our formulation turns out to be quite robust from this point of view since, in some situations, optimal results can be achieved even if the coupling term is computed in approximate way. Moreover, the presence of small cut cells does not affect the stability and conditioning. \cf Next, we discuss the design of effective block preconditioners, whose action is performed thorough direct inversion~\cite{boffi2022parallel} or by employing multigrid algorithms~\cite{mu-interface}. Another issue is related to the fact that the solution of the interface problem can present singularities along the interface. A way to achieve optimal convergence properties is to employ adaptive mesh refinement based on suitable error indicators. This crucial aspect has been previously analyzed in e.g.~\cite{adaptive-1,adaptive-2,adaptive-3,adaptive-4} for other numerical schemes. Thus, we report the \textit{a posteriori} error analysis presented in~\cite{najwa-posteriori} for our FD-DLM formulation.
	
	The second section of the present paper extends the FD-DLM formulation and the related techniques to fluid-structure interaction problems, where an incompressible structure made of viscous hyper-elastic material is immersed in an incompressible Newtonian fluid (see~\cite{2015} and references therein). In this case, the regions occupied by the fluid and the solid evolve in time and their configuration is itself an unknown of the problem. The actual position of the solid is obtained by mapping the reference domain through the deformation unknown: this gives further difficulties for the computation of the so-called coupling terms. Viscous hyper-elastic materials are useful when dealing with several applications: we mention, for instance, blood flow in heart~\cite{peskin1972flow,peskin1989three} and vessels~\cite{delfino1997residual,holzapfel2000new,balzani2016}, and flapping flexible filaments in soap films~\cite{soap1,soap2}. This choice of solid constitutive law is relevant for our formulation since it allows to separate the viscous and elastic parts of the Cauchy stress tensor so that we can treat each of them in the appropriate Eulerian and Lagrangian frameworks, respectively.
	
	Both sections~\ref{sec:interface} and~\ref{sec:fsi} are completed by some numerical results that confirm the theoretical findings and show the versatility of the proposed method. \fc Moreover, several comments throughout this work are devoted to compare the main properties of our formulation with those of other unfitted techniques. For a complete and systematic comparison, we refer the reader to~\cite{boffi2023comparison}. \cf
	
	\section*{Notation}
	
	Given an open bounded domain $D$, we denote by $L^2(D)$ the space of square integrable functions on $D$, endowed with the norm $\|\cdot\|_{0,D}$ associated with the inner product $(\cdot,\cdot)_D$. The subspace $L^2_0(D)$ includes null mean-valued functions.
	
	Sobolev spaces are denoted by the standard symbol $W^{s,p}(D)$, where $s\in\R$ refers to differentiability, while $p\in[1,\infty]$ is the integrability exponent. If $p=2$, then the notation $H^s(D)=W^{s,2}(D)$ is employed, with associated norm $\|\cdot\|_{s,D}$ and semi-norm $|\cdot|_{s,D}$.
	
	Vector-valued functions and spaces will be indicated by boldface letters.
	
	\section{Fictitious domain formulation for elliptic interface problems}~\label{sec:interface}
	
\begin{figure}
		\begin{minipage}[c]{.25\linewidth}
			\begin{center}
				\includegraphics[width=1\linewidth]{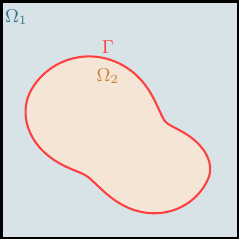} 
			\end{center}
		\end{minipage}\qquad
		\begin{minipage}[c]{.25\linewidth}
			\begin{center}
				\includegraphics[width=1\linewidth]{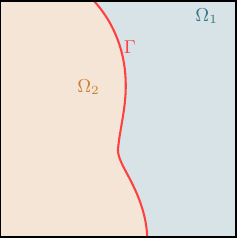}
			\end{center}
		\end{minipage}
		\caption{Two possible domain configurations for the interface problem. Left: $\Omega_2$ is completely immersed in $\Omega$. Right: $\Omega_2$ is not immersed.}
		\label{fig:interface}
	\end{figure}

	Let $\Omega \subset \mathbb{R}^d$, with $d=1,2,3$, be an open bounded domain with Lipschitz boundary~$\partial \Omega$. The domain $\Omega$ is made of two disjoint regions, $\Omega_1$ and $\Omega_2$, separated by the interface ${\Gamma = \overline{\Omega}_1 \cap \overline{\Omega}_2}$. We assume $\Gamma$ to be Lipschitz continuous. Two examples of domain configuration are sketched in Figure~\ref{fig:interface}.
	
	
	We consider the following elliptic interface problem with jumping coefficients.
	\begin{pro}\label{pro:strong_elliptic}
		Given $f_1:\Omega_1\rightarrow\R$, $f_2:\Omega_2\rightarrow\R$, find $u_1:\Omega_1\rightarrow\R$ and $u_2:\Omega_2\rightarrow\R$ satisfying
		\begin{equation*}
			\begin{aligned}
				-\div\left( \nu_1 \grad u_1 \right) &= f_1&&\text{in }\Omega_1\\
				-\div\left( \nu_2 \grad u_2 \right) &= f_2&& \text{in } \Omega_2\\
				u_1-u_2 &= 0&& \text{on } \Gamma \\
				\nu_1 \, \nabla u_1 \cdot \n_1 + \nu_2 \, \nabla u_2 \cdot \n_2 &= 0 &&\text{on } \Gamma\\
				u_1 &= 0&&\text{on } \partial \Omega_1\setminus\Gamma\\
				u_2 &= 0&&\text{on } \partial \Omega_2\setminus\Gamma.
			\end{aligned}
		\end{equation*}
	\end{pro}
	\noindent The notation $\n_i$ ($i=1,2$) refers to the unit normal vector to $\Gamma$, pointing outward from $\Omega_i$. Transmission conditions enforce the continuity of $u_1$ and $u_2$, as well as the continuity of co-normal derivatives, along the interface $\Gamma$. For simplicity,  we assume that the coefficients $\nu_1$ and $\nu_2$ are positive constants.
	
	By setting
	\begin{equation*}
		\tilde\nu = \begin{cases}
			\nu_1&\text{ in }\Omega_1\\
			\nu_2&\text{ in }\Omega_2,
		\end{cases}
		\qquad
		\tilde f = \begin{cases}
			f_1&\text{ in }\Omega_1\\
			f_2&\text{ in }\Omega_2,
		\end{cases}
	\end{equation*}
	the variational formulation of Problem~\ref{pro:strong_elliptic} can be written by applying standard arguments as follows:

	Given $\tilde f\in L^2(\Omega)$, find $\tilde u\in H^1_0(\Omega)$ such that
	\begin{equation}\label{eq:plain_variational}
		\begin{aligned}
		&(\tilde \nu\grad \tilde u,\grad v)_{\Omega}
		= (\tilde f,v)_{\Omega}
		\qquad\forall v \in H^1_0(\Omega).
		\end{aligned}
	\end{equation}
	Equation~\eqref{eq:plain_variational} is equivalent to Problem~\ref{pro:strong_elliptic} with $\tilde u_{|\Omega_i}=u_i$ with $i=1,2$.

	\begin{rem}\label{rem:regularity}
	In general, due to the jump of the normal derivative of $\tilde u$ on $\Gamma$, the solution of~\eqref{eq:plain_variational} does not belong to $H^2(\Omega)$. Indeed, it is well-known that the regularity of the solution to elliptic interface problems with discontinuous coefficients on a domain with Lipschitz boundary belongs to $H^r(\Omega)$ with $1<r<3/2$. However, it is possible to show that the restriction $\tilde u_{|\Omega_i}$, $i=1,2$, is more regular depending on the smoothness of $\Gamma$. If $\Gamma$ is only Lipschitz continuous and presents re-entrant corners, then $u_i \in H^s(\Omega_i)$ with $3/2 < s \le 2$ (see~\cite{nicaise}).
	\end{rem}
	
	Therefore, solving~\eqref{eq:plain_variational} by the finite element method on a mesh which does not take into account the presence of $\Gamma$, gives a non-optimal order of convergence. An optimal method can be recovered by constructing a mesh fitted with the interface $\Gamma$. Clearly, this requirement is strong, especially if the problem under consideration derives from the time discretization of a more complex one, such as a fluid-structure interaction problem. In such case, the interface $\Gamma$ may evolve in time and, thus the mesh should be updated according to the position of $\Gamma$ at each time step. Our approach is to use the fictitious domain formulation with distributed Lagrange multiplier.
	
	\subsection{Fictitious domain formulation}\label{sec:continuous_elliptic}
	In order to reformulate the problem within the fictitious domain framework with distributed Lagrange multiplier (FD-DLM), we extend $\Omega_1$ to the entire $\Omega$, by incorporating the region occupied by $\Omega_2$. For more details, see~\cite{auricchio2015fictitious,boffi2014mixed}.
	
	In detail, $\nu_1$  and $f_1$, originally defined in $\Omega_1$, are prolonged to $\Omega$. We denote the corresponding extensions by $\nu$ and $f\in L^2(\Omega)$, setting $\nu_{|\Omega_1}=\nu_1$ and $f_{|\Omega_1}=f_1$, respectively. 
	
	We denote by $u\in H^1_0(\Omega)$ the extension of $u_1$ to $\Omega$ satisfying $u_{|\Omega_1} = u_1$. The extended $u$ is required to match $u_2\in H^1(\Omega_2)$ in $\Omega_2$. We enforce this superimposition between $\urd$ and $u_2$ at variational level by introducing a suitable functional space $\Lambda$ and a bilinear form $c:\Lambda\times H^1(\Omega_2)\rightarrow\R$ satisfying
	\begin{equation}\label{eq:coupling_def}
		c(\mu,v_2) = 0\quad\forall \mu\in\Lambda
		\qquad
		\text{implies}
		\qquad
		v_2 = 0\quad\text{in }\Omega_2.
	\end{equation}
	We are considering two possible definitions of $\Lambda$ and $c(\cdot,\cdot)$:
	\begin{itemize}
		\item $\Lambda=\Vtd$, the dual space of $\Vt$, and $c(\mu,v_2)=\langle \mu, v_2\rangle$ the duality pairing between $\Vt$ and $\Vtd$;
		\item $\Lambda=\Vt$ and $c(\mu,v_2)=(\mu,v_2)_{\Omega_2}+(\grad\mu,\grad v_2)_{\Omega_2}$ the associated inner product.
	\end{itemize}
	 The two options provide equivalent formulations at the continuous level, since for any $\mu\in\Lambda$, there exists $\varphi\in H^1(\Omega_2)$ being the solution of
	 $$
	 \int_{\Omega_2} (\varphi\cdot v_2 + \grad\varphi\cdot\grad v_2)\, dx_2 = \langle \mu, v_2 \rangle\qquad\forall v_2\in \Vt
	 $$
	 with $\|\varphi\|_{1,\Omega_2} = \|\mu\|_{\Lambda}$.
	
	We introduce the Lagrange multiplier $\lambda\in\Lambda$ associated with the condition $u_{|\Omega_2}=u_2$, then the weak formulation of Problem~\ref{pro:strong_elliptic} in a fictitious domain setting is:
	\begin{pro}\label{pro:fd_ell}
		Given $f\in\LdO$, $f_2\in\LdOdd$, $\nu$, and $\nu_2$, find $(u,u_2,\lambda)\in \Huzo \times \Vt \times \Lambda$ such that
		\begin{equation*}
			\begin{aligned}
				&(\nu\grad u,\grad v)_\Omega + c(\lambda,\vrd) = (f,v)_\Omega&&\forall v\in \Huzo\\
				& ((\nu_2-\nu)\grad u_2,\grad v_2)_{\Omega_2} - c(\lambda,v_2)=(f_2-f,v_2)_{\Omega_2}&&\forall v_2\in \Vt\\
				& c(\mu,\urd-u_2)=0&&\forall\mu\in\Lambda.
			\end{aligned}
		\end{equation*}
	\end{pro}
	
	In~\cite[Thm. 2]{auricchio2015fictitious}, it has been proved that Problem~\ref{pro:fd_ell} is equivalent to the standard variational formulation in~\eqref{eq:plain_variational}.
	The following proposition states the well-posedness of Problem~\ref{pro:fd_ell}, see~\cite[Prop. 1]{auricchio2015fictitious}.
	\begin{prop} \label{prop:well-posed}
		Let $\nu$ and $\nu_2$ be two positive constants. Given $f\in\LdO$ and $f_2\in\LdOdd$, Problem~\ref{pro:fd_ell} admits a unique solution ${(u,u_2,\lambda)\in\Huzo\times\Vt\times\Lambda}$ satisfying
		$$
		\norm{u}_{1,\Omega}+\norm{u_2}_{1,\Omega_2}+\norm{\lambda}_{\Lambda} \leq C \left( \norm{f}_{0,\Omega}+\norm{f_2}_{0,\Omega_2}\right).
		$$
	\end{prop}

	Since the problem shows a saddle point structure, its well-posedness is established by proving the following sufficient conditions~\cite{mixedFEM}:
	\begin{itemize}
		\item \textbf{Ellipticity in the kernel:} there exists a constant $\zeta_1>0$ such that
		\begin{equation}\label{eq:elker}
		(\nu\grad u,\grad v)_\Omega +  ((\nu_2-\nu)\grad u_2,\grad v_2)_{\Omega_2} \ge \zeta_1\,(\|v\|_{1,\Omega}^2 + \|v_2\|_{1,\Omega_2}^2)
		\qquad \forall (v,v_2)\in\mathbb{K}_{\Cmatr},
		\end{equation}
		where $\mathbb{K}_{\Cmatr} = \{(v,v_2)\in \Huzo\times\Vt:\,c(\mu,\vrd-v_2)=0\ \forall\mu\in\Lambda\}$;
		\item \textbf{Inf-sup:} there exists a constant $\zeta_2>0$ such that
		\begin{equation}\label{eq:infsup}
		\sup_{(v,v_2)\in \Huzo\times\Vt} \frac{c(\mu,\vrd-v_2)}{\sqrt{\|v\|_{1,\Omega}^2 + \|v_2\|_{1,\Omega_2}^2}} \ge \zeta_2\,\|\mu\|_\Lambda
		\qquad \forall\mu\in\Lambda.
		\end{equation}
	\end{itemize}

	\begin{rem} \label{rem:stability_c_FDDLM}
		The stability estimate in Proposition~\ref{prop:well-posed} depends on the choice of~$\Lambda$ and~$c(\cdot,\cdot)$ through the constant $C$.
	\end{rem}

	\begin{rem}
		Notice that, due to the equivalence of Problem~\ref{pro:fd_ell} and~\eqref{eq:plain_variational}, we have that the extended solution $u$ belongs to $H^r(\Omega)$ with $1<r<3/2$. On the other hand, $u_2$ belongs to $H^s(\Omega_2)$ with $3/2<s\le2$, see Remark~\ref{rem:regularity}.
	\end{rem}

	\subsection{Discretization} \label{sec:discrete_elliptic}
	We discretize Problem~\ref{pro:fd_ell} by mixed finite elements. We partition $\Omega$ and $\Omega_2$ by two independent shape regular meshes $\T_h$ and $\T_h^2$, respectively. We denote by $h$ the size of $\T_h$, whereas $h_2$ denotes the mesh size of $\T_h^2$. The arguments we are going to present in this section are valid for both triangular and quadrilateral meshes, \red although we present in detail the quadrilateral case.\nored

	We define the discrete finite element spaces as follows: $V_h \subset \Huzo$, $\Vth \subset \Vt$, and $\Lambda_h \subset \Lambda$. At the discrete level, if $\Lambda=\Vtd$, the duality pairing between $\Vt$ and its dual can be evaluated using the scalar product in $\LdOdd$, provided that $\Lambda_h\subset\LdOdd$. Thus we have
	\begin{equation}\label{eq:coupling_l2}
		c(\mu_h,\vth) = (\mu_h,\vth)_{\Omega_2}\qquad\forall\mu_h\in\Lambda_h,\,\vth\in \Vth.
	\end{equation}
	On the other hand, if $\Lambda=\Vt$, the inner product in $\HuOd$ can be used also at discrete level if $\Lambda_h$ contains continuous functions.
	
	The discrete formulation of Problem~\ref{pro:fd_ell} reads as follows.
	\begin{pro}\label{pro:fd_ell2}
		Given $f\in\LdO$, $f_2\in\LdOdd$, $\nu$, and $\nu_2$, find $(u_h,\uth,\lambda)\in V_h \times \Vth \times \Lambda_h$ such that 
		\begin{equation*}
			\begin{aligned}
				&(\nu\grad u_h,\grad v_h)_\Omega + c(\lambda_h,\vhrd) = (f,v_h)_\Omega&&\forall v_h\in V_h\\
				& ((\nu_2-\nu)\grad \uth,\grad \vth)_{\Omega_2} - c(\lambda,\vth)=(f_2-f,\vth)_{\Omega_2}&&\forall \vth\in \Vth\\
				& c(\mu_h,\uhrd- \uth)=0&&\forall\mu_h\in\Lambda_h.
			\end{aligned}
		\end{equation*}
	\end{pro}

	For a generic element $\E$, we denote by $\Poly_1(\E)$ the space of linear polynomials if $\E$ is a simplex and by $\Qoly_1(\E)$ the space of bilinear polynomials if $\E$ is a quadrilateral or an hexahedron. Moreover, let $\Poly_0(\E)$ be the space of constants and $\Boly(\E)$ the space of bubbles vanishing at the boundary of an element $\E$. 
	
	We solve Problem~\ref{pro:fd_ell2} by considering two distinct sets of finite element spaces, which are distinguished by the Lagrange multiplier space being continuous or discontinuous depending on the choice of the coupling term. More precisely, we consider continuous linear elements for all variables~\cite{auricchio2015fictitious}, that is
	\begin{subequations}\label{eq:elm1}
		\begin{align}
			&V_h  = \left\{ v_h \in V : v_{h|\E} \in \Qoly_1(\E), \forall \E \in \T_h \right\},\\
			&\Vth = \left\{ \vth \in V_2 : v_{2,h|\Ed} \in \Qoly_1(\Ed), \forall \Ed \in \T_h^2 \right\},\\
			&\Lambda_h = \Vth
		\end{align}
	\end{subequations}
	or, if the coupling term is defined as~\eqref{eq:coupling_l2}, we take
	\begin{subequations}\label{eq:elm2}
	\begin{align}
		&V_h = \left\{ v_h \in V: v_{h|\E} \in \Qoly_1(\E), \forall \E \in \T_h \right\} \subset V, \\
		&\Vth = \left\{ \vth \in V_2 : v_{2,h|\Ed} \in \Qoly_1(\Ed) \oplus \Boly(\Ed), \forall \Ed \in \T_h^2 \right\} \subset V_2, \\
		&\Lambda_h = \left\{ \lambda_h \in \Lambda : \lambda_{h|\Ed} \in \Poly_0(\Ed), \forall \Ed \in \T_h^2 \right\} \subset \Lambda.
	\end{align}
	\end{subequations}
	For the second choice the space $\Vth$ is enriched by the local bubble space $\Boly$ to ensure the stability of the discrete formulation. Bi-quadratic elements for $V_h$ and $\Vth$ could be considered as well~\cite{najwa}.
	
	Both choices are stable since they satisfy the discrete counterpart of~\eqref{eq:elker} and~\eqref{eq:infsup}. In particular, stability of~\eqref{eq:elm1} has been proved for both simplicial and quadrilateral/hexahedral meshes in~\cite{auricchio2015fictitious,boffi2014mixed}. The ellipticity in the kernel condition is met under the assumption that $\nu_2 > \nu$. This requirement can be relaxed for simplicial meshes if $h_2/h^d$ is small enough, see~\cite{boffi2014mixed}. Stability of~\eqref{eq:elm2} has been proved in~\cite{najwa}. The ellipticity in the kernel condition is satisfied assuming that $\nu_2 > \nu$. Numerical findings suggest that this requirement might be relaxed~\cite{najwa}.

	The error estimates on the discrete solution are provided in~\cite{auricchio2015fictitious,boffi2014mixed,najwa}. We point out that such bounds depend on the choice of~$\Lambda_h$. The following estimate holds.
	\begin{prop}\label{prop:error_estimate}
		Let $(u,u_2,\lambda)$ and $(u_h,u_{2h},\lambda_h)$ be the solutions of Problems~\ref{pro:fd_ell} and~\ref{pro:fd_ell2}, respectively. If the discrete spaces are given by~\eqref{eq:elm1}, then the following error estimate holds true
		\begin{align*}
			\norm{u-u_h}_{1,\Omega}&+\norm{u_2-\uth}_{1,\Omega_2}+\norm{\lambda-\lambda_h}_{\Lambda}\\
			&\leq C_1\left(  h^{r-1}\norm{u}_{r,\Omega} +h_2^{s-1}\norm{u_2}_{s,\Omega_2} + h_2^{s-1}
			\norm{({\nu}/{\nu_2})f_2-f}_{0,\Omega_2} \right),
		\end{align*}
		while, in the case of the spaces in~\eqref{eq:elm2}, it holds
		\begin{align*}
			\norm{u-u_h}_{1,\Omega}&+\norm{u_2-\uth}_{1,\Omega_2}+\norm{\lambda-\lambda_h}_{\Lambda}\\
			&\leq C_2\left(  h^{r-1}\norm{u}_{r,\Omega} +\max(h_2^{s-1},h_2^{1-t})\norm{u_2}_{s,\Omega_2} + h_2
			\norm{({\nu}/{\nu_2})f_2-f}_{0,\Omega_2} \right),
		\end{align*}
		where $r$, $s$ are given in Remark~\ref{rem:regularity}, $1/2<t<1$ and $C_1,\,C_2$ are positive constants independent of $h$ and~$h_2$.
	\end{prop}
	

	\subsection{How to deal with the coupling term} \label{sec:coupling}
	As we previously mentioned, the main idea behind the fictitious domain consists in extending the domain $\Omega_1$ to the region occupied by $\Omega_2$. The domain~$\Omega_2$ is then superimposed to the background $\Omega$. Such superimposition is mathematically enforced by means of the Lagrange multiplier $\lambda$, which is applied through the coupling bilinear form~$c(\cdot,\cdot)$. We notice that, at continuous level, the bilinear form $c(\cdot,\cdot)$ has the same representation when applied to functions in $\Vt$ and to the restriction of elements in $\Huzo$. On the other hand, at discrete level, Problem~\ref{pro:fd_ell2} in matrix form reads
	\begin{equation}\renewcommand{\arraystretch}{1.5}
		\label{eq:matrix}
		\left[\begin{array}{@{}cc|c@{}}
			\Amatr & \Zmatr & \Cmatr_1^\top \\
			\Zmatr & \Amatr_2 & -\Cmatr_2^\top \\
			\hline
			\Cmatr_1 & -\Cmatr_2 & \Zmatr
		\end{array}\right]
		\left[ \begin{array}{c}
			u\\
			u_2\\
			\hline
			\lambda
		\end{array}\right]=
		\left[ \begin{array}{c}
			\mathsf{f}\\
			\mathsf{f}_2-\mathsf{f}\\
			\hline
			0
		\end{array}\right].
	\end{equation}\renewcommand{\arraystretch}{1.15}
	Let us denote by $\varphi$, $\psi$ and $\zeta$ the basis functions of $V_h$, $\Vth$ and $\Lambda_h$, respectively. We have
	\begin{equation}
		\begin{aligned}
			&(\Amatr)_{i,j} = (\nu\grad \varphi_j,\grad \varphi_i)_\Omega,\quad
			&&(\Amatr_2)_{i,j} = ((\nu_2-\nu)\grad \psi_j,\grad \psi_i)_{\Omega_2},\\
			& (\Cmatr_1)_{i,k} = c(\zeta_k,(\varphi_i)_{|\Omega_2}),\quad
			&& (\Cmatr_2)_{i,k} = c(\zeta_k,\psi_i).
		\end{aligned}
	\end{equation}

	We highlight that the matrices associated with the bilinear form $c(\cdot,\cdot)$ are denoted differently because they are constructed in two different ways.
	
\begin{figure}
		\centering
		\subfloat[\label{fig:mismatch}]{\includegraphics[width=0.305\linewidth]{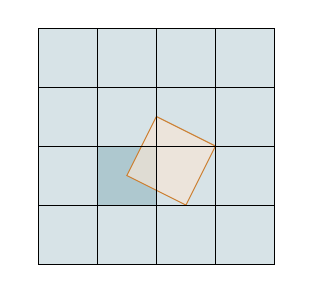}}\quad
		\subfloat[\label{fig:intersection}]{\includegraphics[width=0.295\linewidth]{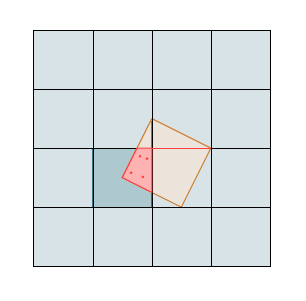}}\quad
		\subfloat[\label{fig:inexact}]{\includegraphics[width=0.3\linewidth]{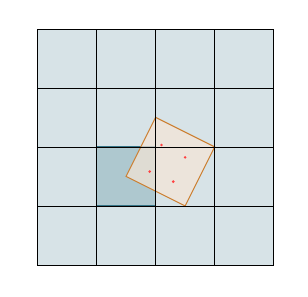}}
		\caption{(\textsc{A}) The immersion of $\Ed\in\T_h^2$ into the background mesh $\T_h$. The immersed element, depicted in beige, is not matching the support of the finite element functions defined on $\T_h$ (blue). (\textsc{B}) The exact computation of the interface matrix requires the implementation of a composite quadrature rule on the intersection $\T_h^2\cap\T_h$. (\textsc{C}) Quadrature points of an inexact rule in $\Ed$.}
		\label{fig:coupling_term}
	\end{figure}

	Since the mesh $\T_h^2$ of $\Omega_2$ is used to discretize both $u_2$ and $\lambda$, the matrix $\Cmatr_2$ is assembled with the standard procedure of finite elements. Hence, we focus on the interface matrix~$\Cmatr_1$. The assembly procedure for this finite element matrix is not trivial since integration over non-matching grids is required (see~\cite{boffi2022interface}). Indeed, the definition of~$\Cmatr_1$ associated with the scalar product either in~$\LdOdd$ or~$\HuOd$ requires the integration over~$\Omega_2$ of~$\zeta_k\in\Lambda_h$, which is defined on~$\Omega_2$, multiplied by~$\varphi_i$, which is a function defined on the background mesh $\T_h$ and then restricted to the immersed domain. An example of configuration is sketched in Figure~\ref{fig:mismatch}: it is clear that the immersed beige element is not matching with the support of the finite element basis functions defined on the background mesh.
	
	Here we describe in detail how to deal with the assembly of $\Cmatr_1$. The coupling term can be written as the sum of local contributions: indeed, if $c(\cdot,\cdot)$ is the scalar product in $\LdOdd$, we have
	\begin{equation}\label{eq:interf_l2}
		c(\mu_h,\vhrd) = \int_{\Omega_2} \mu_h\cdot\vhrd\,dx_2 = \sum_{\Ed\in\T_h^2} \int_{\Ed} \mu_h\cdot\vhrd\,dx_2,
	\end{equation}
	while, if $c(\cdot,\cdot)$ is the scalar product in $\HuOd$, it holds
	\begin{equation}\label{eq:interf_h1}
		\begin{aligned}
		c(\mu_h,\vhrd) &= \int_{\Omega_2} \mu_h\cdot\vhrd + \grad\mu_h\cdot\grad\vhrd \,dx_2\\
		&= \sum_{\Ed\in\T_h^2} \left[\int_{\Ed}\mu_h\cdot\vhrd\,dx_2 + \int_{\Ed}\grad\mu_h\cdot\grad\vhrd\,dx_2\right].
		\end{aligned}
	\end{equation}
	We are going to discuss how to compute the local integrals reported in the above equations. We will focus on two possible techniques. Indeed, the exact computation of such quantities requires the design of a composite quadrature rule on the intersection between $\Ed\in\T_h^2$ and the background mesh~$\T_h$. We could also use a quadrature rule directly on each $\Ed\in\T_h^2$, which leads to a quadrature error. In the remainder of this section we restrict the discussion to the two-dimensional case in order to keep the presentation simpler.
	
	\subsubsection{Exact procedure}
	
	We first describe the exact procedure, which requires a composite quadrature rule on~$\Ed$. By computing the intersection between the considered immersed element and the background mesh, we obtain a partition of $\Ed$ into $J$ disjoint polygons
	\begin{equation*}
		\Ed = P_1 \cup \cdots \cup P_J,\qquad\qquad
		\text{with}\quad P_{j_1}\cap P_{j_2}=\emptyset\quad\text{ if }j_1\neq j_2,
	\end{equation*}
	as depicted in Figure~\ref{fig:intersection}.
	More precisely, each $P_j$ is contained in a single element $\E\in\T_h$, so that the background basis functions are completely supported in $P_j$. We have at least two options for integrating in each $P_j$. If $P_j$ is not a triangle or a quadrilateral, we construct an auxiliary triangulation $\{T_i\}_{i=1,\dots,I(j)}$ by simply connecting the barycenter with the vertices. This gives
	\begin{equation*}
		\begin{aligned}
			& \int_{\Ed} \mu_h\cdot\vhrd\,dx_2
			= \sum_{j=1}^J \int_{P_j} \mu_h\cdot\vhrd\,dx_2\\
			&\hspace{2.7cm}= \sum_{j=1}^J \sum_{i=1}^{I(j)} \int_{T_i} \mu_h\cdot\vhrd\,dx_2
			= \sum_{j=1}^J \sum_{i=1}^{I(j)}\mathrm{area}(T_i)\sum_{\ell=1}^{L_0} \quadweigth_\ell^0\,\mu_h(\quadnode_\ell^0)\cdot\vhrd(\quadnode_\ell^0)
		\end{aligned}
	\end{equation*}
	and
	\begin{equation*}
		\begin{aligned}
			& \int_{\Ed} \grad\mu_h\cdot\grad\vhrd\,dx_2
			= \sum_{j=1}^J \int_{P_j} \grad\mu_h\cdot\grad\vhrd\,dx_2= \sum_{j=1}^J \sum_{i=1}^{I(j)} \int_{T_i} \grad\mu_h\cdot\grad\vhrd\,dx_2\\
			&\hspace{3.4cm}
			= \sum_{j=1}^J \sum_{i=1}^{I(j)} \mathrm{area}(T_i)\sum_{\ell=1}^{L_1} \quadweigth_\ell^1\,\grad\mu_h(\quadnode_\ell^1)\cdot\grad\vhrd(\quadnode_\ell^1),
		\end{aligned}
	\end{equation*}
	where we denoted by $\{(\quadnode_\ell^0,\quadweigth_\ell^0)\}_{\ell=1,\dots,L_0}$ and $\{(\quadnode_\ell^1,\quadweigth_\ell^1)\}_{\ell=1,\dots,L_1}$ nodes and weights of two not necessarily equal quadrature rules.
	
	Alternatively, polygonal quadrature rules may be directly employed to avoid the construction of the auxiliary triangulation, see~\cite{polygauss} for instance. In any case, all quantities are computed exactly provided that a precise enough quadrature rule is employed either in each $T_i$ ($i=1,\dots,I(j)$) or in~$P_j$ ($j=1,\dots,J$). 
	
	\subsubsection{Inexact procedure}
	
	The computation of the intersection between the involved meshes may be expensive, unless special computational techniques are employed for an efficient implementation.
	
	A cheaper procedure can be then employed for computing the quantities in~\eqref{eq:interf_l2} and~\eqref{eq:interf_h1}. Indeed, we can skip the computation of the mesh intersection and integrate in each $\Ed\in\T_h^2$ using a certain quadrature rule, i.e.
	\begin{equation*}
		\begin{aligned}
			& \int_{\Ed} \mu_h\cdot\vhrd\,dx_2
			\approx \mathrm{area}(\Ed)\sum_{\ell=1}^{L_0} \quadweigth_\ell^0\,\mu_h(\quadnode_\ell^0)\cdot\vhrd(\quadnode_\ell^0)\\
			& \int_{\Ed} \grad\mu_h\cdot\grad\vhrd\,dx_2
			\approx \mathrm{area}(\Ed)\sum_{\ell=1}^{L_1} \quadweigth_\ell^1\,\grad\mu_h(\quadnode_\ell^1)\cdot\grad\vhrd(\quadnode_\ell^1)\\
		\end{aligned}
	\end{equation*}
	
	An example is sketched in Figure~\ref{fig:inexact}, where we highlight four quadrature nodes in the immersed element. Hence, the nodes are placed in different elements of the background mesh $\T_h$. When evaluating the shape functions defined on $\T_h$, we first check the location of the quadrature points and then we evaluate accordingly.
	
	Since we are integrating piecewise polynomials, this procedure is inexact and originates a quadrature error.
	
	\begin{rem}
		In~\cite{BCG24}, we studied the behavior of the quadrature error in the case of FSI problems discretized by piecewise linear elements on triangular meshes, however the same arguments are still valid in the case of simpler elliptic interface problems. We proved that if $c(\cdot,\cdot)$ is the scalar product in $\LdOdd$, then the quadrature error is comparable to the other discretization errors. Conversely, if the $\HuOd-$scalar product is considered, then the quadrature error is large unless the ration $h_2/h$ is decaying fast enough. The interested reader may refer to~\cite{BCG24} for a complete discussion about these aspects. 
	\end{rem}

	\begin{rem}
		We discussed the assembly techniques for the interface matrix in the case of two dimensional problems, where triangular or quadrilateral meshes are employed. The same computational techniques can be considered when solving three dimensional problems. We have not proved convergence estimates for the quadrature error in the case of 3D problems, but we expect the same behavior described in the previous remark. 
	\end{rem}

	\begin{rem}
		It is important to observe that the efficient integration of the coupling terms cannot rely on standard algorithms~\cite{krause2016parallel,boffi2022interface,boffi2023comparison}. 
		For instance, when computing the intersection between $\T_h$ and $\T_h^2$, the total number of possible matching is given by $\mathcal{N}(\T_h)\times\mathcal{N}(\T_h^2)$, where $\mathcal{N}(\T_h)$ and $\mathcal{N}(\T_h^2)$ denote the number of elements of $\T_h$ and $\T_h^2$, respectively. A naive implementation of such procedure would be prohibitive when dealing with fine meshes due to the large amount of required testing operations.
		In order to reduce the computational cost, a \textit{collision detection algorithm}~\cite{collision} can be employed to efficiently identify pairs of intersecting cells or background cells where quadrature nodes may be located. A popular technique is based on bounding boxes within the framework of \verb*|R-tree| data structures, where each mesh can be interpreted as a hierarchy of objects~\cite{rtree1,rtree2}. \verb*|R-tree| structures support the efficient resolution of boolean and nearest-neighbor geometric operations.
	\end{rem}

	\begin{rem}
		It is well-known that unfitted finite elements may require additional terms to stabilize the discrete formulation in the presence of small cut cells: this is the case, for instance, of Cut--FEM~\cite{burman1,burman2,wadbro2013uniformly,HANSBO201490} and Finite Cells~\cite{finitecell1,finitecell2,finitecell3} methods. We point out that our formulation is naturally stable without the need of adding any artificial penalization term~\cite{stat,BCG25}.
	\end{rem}
		
	\subsection{Solvers}\label{sec:solvers}
	
		The linear system~\eqref{eq:matrix} arising from Problem~\ref{pro:fd_ell2} is usually large and its direct solution is not feasible. The use of an iterative solver is thus recommended and GMRES represents a good option, although the design of effective preconditioners is not straightforward. 
		We presented and discussed a first parallel solver for the fictitious domain formulation of fluid-structure interaction problems in~\cite{mca28020059,boffi2022parallel}, where we performed an in depth numerical investigation to assess the robustness of the solver with respect to mesh refinement and strong scalability. In these works, we proposed a diagonal and a triangular preconditioner showing that the triangular one is better performing.
		
		The preconditioners introduced in~\cite{mca28020059,boffi2022parallel} can also be considered for the solution of the linear system~\eqref{eq:matrix}. By separating the variables defined in $\Omega$ and $\Omega_2$, respectively, the global matrix can be subdivided into four blocks as
		\begin{equation*}\renewcommand{\arraystretch}{1.5}
			\left[\begin{array}{@{}c|cc@{}}
				\Amatr & \Zmatr &
				 \Cmatr_1^\top \\
				 \hline
				\Zmatr & \Amatr_2 & -\Cmatr_2^\top \\
				\Cmatr_1 & -\Cmatr_2 & \Zmatr
			\end{array}\right],
		\end{equation*}\renewcommand{\arraystretch}{1.15}
		so that the following preconditioners can be used
		\begin{equation}\label{eq:prec}\renewcommand{\arraystretch}{1.5}
			\mathsf{P}_{\mathsf{diag}} = \left[\begin{array}{@{}c|cc@{}}
				\Amatr & \Zmatr &
				\Zmatr \\
				\hline
				\Zmatr & \Amatr_2 & -\Cmatr_2^\top \\
				\Zmatr & -\Cmatr_2 & \Zmatr
			\end{array}\right],
			\qquad\qquad
			\mathsf{P}_{\mathsf{tri}} = \left[\begin{array}{@{}c|cc@{}}
				\Amatr & \Zmatr &
				\Zmatr \\
				\hline
				\Zmatr & \Amatr_2 & -\Cmatr_2^\top \\
				\Cmatr_1 & -\Cmatr_2 & \Zmatr
			\end{array}\right].\renewcommand{\arraystretch}{1.15}
		\end{equation}
	Examining the structure of the blocks in $\mathsf{P}{\mathsf{diag}}$ and $\mathsf{P}{\mathsf{tri}}$, we observe that the matrix inversion depends on the invertibility of the blocks $\Amatr$ and $\mathsf{L}=[\Amatr_2,-\Cmatr_2^\top;-\Cmatr_2,\Zmatr]$. As shown in~\cite{alshehri2025multigrid}, these blocks are invertible for both finite element discretization choices and for both $\Lambda$ space options.
	In our first implementation, $\mathsf{P}_{\mathsf{diag}}$ and $\mathsf{P}_{\mathsf{tri}}$ were inverted by a direct solver, but for large problems this may not be feasible and iterative algorithms should be preferred. The behavior of the block $\Amatr$ does not give any particular issue, whereas the lower right block $\mathsf{L}$ is more challenging. 	
	%
	%
	
	Inexact inversion of such preconditioners using multigrid methods is discussed in~\cite{alshehri2025multigrid}. The matrix $\Amatr$ is symmetric positive definite and can be inverted using multigrid methods with standard smoothers such as Jacobi or SOR methods. The block $\mathsf{L}$, however, is a saddle point matrix that is also invertible, but due to the presence of a zero diagonal block, common smoothing schemes like SOR cannot be used because they involve division by the diagonal entries, which are zero. More advanced strategies, such as a Vanka smoother, should be considered for this purpose.
	The multigrid method is particularly effective in reducing the number of GMRES iterations required for convergence, especially when the block $\Amatr$ is inverted inexactly. This approach notably improves the performance of the algorithm even when $\mathsf{P}_{\mathsf{diag}}$ is used.
	
	\subsection{A posteriori error analysis}
	
	As discussed in Remark~\ref{rem:regularity}, the interface problem under consideration is characterized by solutions with reduced regularity. While quasi-optimal convergence rates can be achieved with uniform mesh refinement, mesh adaptivity remains essential for recovering optimal rates by concentrating computational effort near singularities and interfaces. 
	
	In this context, residual-based estimators have been thoroughly investigated (see \cite{najwa-posteriori}), demonstrating both reliability and efficiency for $c(\cdot,\cdot)$ being the scalar product in $\LdOdd$. 
	
		To properly account for the variability in the source data $f \in L^2(\Omega)$ and $f_2 \in L^2(\Omega_2)$, we introduce oscillation terms that capture the mismatch between the continuous data and their numerical approximations. Let $\Pi_0$ and $\Pi_0^2$ be the $L^2$-projection operators mapping functions onto the space $\Poly_0(\T_h)$ and $\Poly_0(\T_h^2)$, respectively. Here, $\Poly_0(\mathfrak{T}_h)$ stands for the space of piecewise constants on the generic mesh $\mathfrak{T}_h$. For $\E\in\T_h$ and $\Ed\in\T_h^2$, we set:
		$$
		osc_{\E} = h_{\E} \norm{f-\Pi_0f}_{0,\E},
		\qquad\qquad
		osc_{\Ed} = h_{\Ed} \norm{((f_2 -f)-\Pi_0^2(f_2 -f))}_{0,\Ed}.
		$$
		
	\subsubsection{Residuals and error indicators}
	We introduce four residuals: $r_{\E}$, $r_{\Ed}$ for elements, and $r_{\edge}$, $r_{\edged}$ for edges. In addition, we introduce $\tilde{\lambda}$ and $\tilde{\lambda}_h$ to represent suitable extensions of the continuous and discrete Lagrange multipliers, see~\cite[Def. 3.1, Def. 3.2]{najwa-posteriori}. More precisely, we extend by zero $\lambda_h$ to $\Omega$. Hence, we have $\tilde{\lambda}_h:\Omega\rightarrow\R$ with $\tilde{\lambda}_{h|\Omega_2}=\lambda_h$ and $\tilde{\lambda}_{h|\Omega\setminus\Omega_2}=0$. The definition of the extension of $\lambda$ is given by $\tilde{\lambda}\in H^{-1}(\Omega)$ such that
	\begin{equation}\label{eq:lambda-tilde}
		\langle \tilde{\lambda},v \rangle_{\Huzo} = \langle \lambda,v_{|\Omega_2} \rangle_{\Vt}.
	\end{equation}
	Here $\langle\cdot,\cdot\rangle_{W}$ denotes the duality pairing between $W$ and its dual space $W^\prime$. We observe that the first equation in Problem~\ref{pro:fd_ell} reads
	\begin{equation*}
		(\nu\grad u,\grad v)_\Omega +  \langle \tilde{\lambda},v \rangle_{\Huzo}
		= (f,v)_\Omega\qquad\forall v\in\Huzo.
	\end{equation*}
	%

	For each element $\E \in \T_h$, the residual is given by
	\begin{equation}
	r_{\E}(u_h, \tilde{\lambda}_h) = \nabla \cdot (\nu \nabla u_h) - \tilde{\lambda}_h + \Pi_0 f.
	\end{equation}
	Similarly, for any $\Ed \in \T_h^2$, we define
	\begin{equation}
	r_{\Ed}(\uth, \lambda_h) = \nabla \cdot \left( (\nu_2 - \nu) \nabla \uth \right) + \lambda_h + \Pi_0^2 (f_2 - f).
	\end{equation}
	
	Let $\mathcal{E}$ and $\mathcal{E}_2$ be the set of edges of $\T_h$ and $\T_h^2$, respectively. For $\edge \in \mathcal{E}$, the residual depends on whether $\edge$ lies in the interior of the mesh or on the boundary $\partial \Omega$:
	\begin{equation}
	r_{\edge}(u_h) =
	\begin{cases} 
	- \jump{ \nu \nabla u_h \cdot \n_1 }_\edge, & \text{if } \edge \text{ is an interior edge},\\
	0, & \text{if } \edge \text{ is on }\partial \Omega.
	\end{cases}
	\end{equation}
	For $\edged \in \mathcal{E}_2$, the definition changes depending on the location of $\edged$
	\begin{equation}
	r_{\edged}(\uth) =
	\begin{cases} 
	- \jump{ (\nu_2 - \nu) \nabla \uth \cdot \n_2 }_{\edged}, & \text{if } \edged \text{ lies in the interior of } \T_h^2, \\
	- (\nu_2 - \nu) \nabla \uth \cdot \n_2, & \text{if } \edged \text{ is on } \Gamma.
	\end{cases}
	\end{equation}
	In these definitions, the jump operator $\jump{\varphi}_e$ represents the jump of $\varphi$ across the edge $e$. 

	The error indicators are constructed by combining weighted norms of the residuals with the oscillation terms, yielding a comprehensive measure of the local error distribution. Then, for each element $\Ed \in \T_h$ and $\Ed \in \T^2_h$, the local error indicators are defined by:
	\begin{align*}
		& \eta_{\E}^2 = h_{\E}^2 \left\| r_{\E}(u_h, \tilde{\lambda}_h) \right\|_{0,\E}^2  
		+ \frac{1}{2} \sum_{\substack{\edge \subset \partial \E \\ \edge \not\subset \partial \Omega}} 
		h_{\edge} \left\| r_{\edge}(u_h) \right\|_{0, \edge}^2,\\[1.5ex]
		&\eta_{\Ed}^2 = h_{\Ed}^2 \left\| r_{\Ed}(\uth, \lambda_h) \right\|_{0,\Ed}^2  
		+ \left\| u_{h|\Omega_2} - \uth \right\|_{1, \Ed}^2\\
		&\qquad\quad + \frac{1}{2} \sum_{\substack{\edged \subset \partial \Ed \\ \edged \not\subset \Gamma}}  
		h_{\edged} \left\| r_{\edged}(\uth) \right\|_{0, \edged}^2  
		+ \sum_{\substack{\edged \subset \partial \Ed \\ \edged \subset \Gamma}}  
		h_{\edged} \left\| r_{\edged}(\uth) \right\|_{0, \edged}^2.
	\end{align*}
	
	The global error estimator is then defined by the aggregation of all local contributions:
	\begin{equation}\label{eq:eta}
		\eta^2 = \sum_{\E \in \T_h} \eta_{\E}^2,
		\qquad
		\qquad
		\eta_2^2 = \sum_{\Ed \in \T^2_h} \eta_{\Ed}^2.
	\end{equation}

	Those indicators are shown to be both reliable and efficient in \cite{najwa-posteriori}, providing a global upper bound of the true error and a local lower bound, respectively. The reliability of the error estimator is established by the following proposition. 

	\begin{prop}
	Let $(u,u_2,\lambda)$ be the exact solution of the continuous problem and $(u_h,\uth,\lambda_h)$ be the solution of the discrete one. Then, there exists a constant ${C} > 0$ that depends only on the domains $\Omega$ and $\Omega_2$, the coefficients $\nu$ and $\nu_2$, and the shape regularity of $\T_h$ and $\T^2_h$, such that
	\begin{equation}
	\label{eq:upper_a_posterori}
	\begin{aligned}
	\|u-u_h\|_{1,\Omega}&+\norm{u_2-\uth}_{1,\Omega_2}+\norm{\lambda-\lambda_h}_{\Lambda}\\
	&\leq {C} \left(
		\sum_{\E \in\T_h} \left( \eta_{\E}^2 + osc_{\E}^2 \right) +
	\sum_{\Ed \in \T^2_h} \left( \eta_{\Ed}^2 + osc_{\Ed}^2 \right) 
	\right) ^{1/2}.
	\end{aligned}
	\end{equation}
	\end{prop}

	The proof, which involve the approximation of functions in $H^1$ by means of a Cl\'ement-type interpolation, exploits the definition of the extension of $\lambda$ given in~\eqref{eq:lambda-tilde}, and that by zero of $\lambda_h$, see \cite[Prop. 3.1]{najwa-posteriori}. This result also holds to the case of non-constant coefficients $\nu$ and $\nu_2$ as described in~\cite[Sect. 3.2]{najwa-posteriori}.

\begin{figure}
		\centering
		\subfloat[$\omega_{\E}$]{\includegraphics[width=0.37\linewidth]{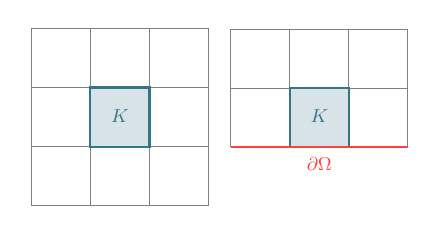}}
		\qquad
		\subfloat[${\omega_{\edge}}$]{\includegraphics[width=0.32\linewidth]{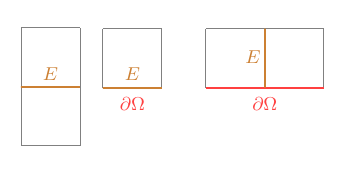}}
		\caption{Examples of element-based $\omega_{\E}$ and edge-based ${\omega_{\edge}}$ neighborhoods}
		\label{fig:neig}
	\end{figure}

	 Before discussing the efficiency of the error estimators, we need to introduce the following neighborhoods of elements and edges as described in Figure~\ref{fig:neig} for a quadrilateral mesh $\T_h$ (similar notation holds for $\T_h^2$):
	\begin{itemize}
		\item \emph{Element-based neighborhoods $\omega_{\E}$ and $\omega_{\Ed}$:} the collection of elements that share a \emph{common edge} with $\E$ and $\Ed$, respectively;
		\item \emph{Edge-based neighborhoods $\omega_{\edge}$ and $\omega_{\edged}$:} the collection of elements for which $\edge$ and $\edged$ belong to, respectively.
	\end{itemize}
	For the efficiency of the error estimators, we have the following result.

	\begin{prop}
		There exist two constants $\underline{C_1}>0$ and $\underline{C_2}>0$, which depend only on the domains $\Omega,~\Omega_2$ and on the shape regularity parameters of $\T_h,~\T^2_h$, such that the local error indicators satisfy the following bounds:
		\begin{equation*}
			\label{eq:LLB}
			\begin{aligned}
				\underline{C}_1\,\eta_{\E}^2 &\leq    
				{\nu}^2\norm{\nabla (u-u_h)}_{0,\omega_{\E}}^2
				+\norm{\tilde{\lambda}-\tilde{\lambda}_h}_{(H^1(\omega_{\E}))^\prime}^2
				+\sum_{\E^\prime \in \omega_{\edge}} osc_{\E^\prime}^2
				,\\
				\underline{C}_2\,\eta_{\Ed}^2&\leq   
				\overline{\nu}^2_2\norm{\nabla (u_2-\uth)}_{0
					,\omega_{\Ed}}^2
				+\norm{\nabla (u-u_h)_{|\Omega_2}}_{0,\omega_{\Ed}}^2
				+\norm{\lambda-\lambda_h}_{(H^{1}(\omega_{\Ed}))^\prime}^2
				+\sum_{\Ed^\prime \in \omega_{\edged}} osc_{\Ed^\prime}^2,
			\end{aligned}
		\end{equation*}
		where $\overline{\nu}_2 = \max\{\nu_2-\nu,1\}$.
	\end{prop}
	
	The proof of this proposition relies on the arguments in~\cite{verfurth2013posteriori}, which make use of bubble functions in order to localize the indicators on elements or edges. We refer to~\cite[Prop. 3.2]{najwa-posteriori} for the details of the proof in the case of constant coefficients $\nu$, $\nu_2$ and \cite[Sect. 3.2]{najwa-posteriori} for the general case.

	\subsection{Adaptive mesh and refinement strategy}\label{sec:adaptive}

Our adaptive refinement is based on the SOLVE--ESTIMATE--MARK--REFINE strategy~\cite{dorfler}. This iterative procedure continues until the global error estimator satisfies the stopping criterion:
\begin{equation*}
    \eta + \eta_2 \leq \texttt{tol},
\end{equation*}
where $\texttt{tol}$ is a user-defined tolerance controlling the accuracy of the numerical solution. The \textit{bulk marking strategy}, or D\"orfler marking~\cite{dorfler}, is used to select elements for refinement. It consists in marking elements whose cumulative error indicators contribute to a fraction $\alpha_1 \in [0, 1]$ of the total error. 

\subsubsection{Adaptive refinement procedure}
We construct a sequence of refinements $\T_{h,\ell}$ and $\T_{h,\ell}^2$ for $\ell=0,\dots$ until the stopping criterion is satisfied. We recall the steps involved in the adaptive algorithm.

\begin{enumerate}[1.]
    \item Initialize meshes $\T_{h,0}$ and $\T^2_{h,0}$ with their respective local error indicators $\eta_{\E,0}$ and $\eta_{\Ed,0}$.
    \item Compute the global error estimators $\eta_0$ and $\eta_{2,0}$ by summing the local contributions.
    \item {While} the stopping criterion $\eta_\ell + \eta_{2,\ell} \le \texttt{tol}$ is \textit{not} satisfied, repeat:
    \begin{enumerate}[a.]
        \item Sort the elements in $\T_{h,\ell}$ and $\T^2_{h,\ell}$ in descending order according to their local error indicators.
        \item Identify a set of elements in $\T_{h,\ell}$ and $\T_{h,\ell}^2$ whose cumulative error indicators account for a fraction $\alpha_1$ of $\eta_\ell$ and $\eta_{2,\ell}$, respectively. Mark these elements for refinement.
        \item Optionally, identify elements with the lowest error indicators contributing to a fraction $\alpha_2$ of $\eta_\ell$ and $\eta_{2,\ell}$, respectively, and mark them for coarsening.
        \item Refine the marked elements by subdividing each $\E\in\T_{h,\ell}$ and $\Ed\in\T^2_{h,\ell}$ into smaller sub-elements.
        \item Coarsen elements marked for coarsening, provided that their neighboring refinement levels remain consistent.
        \item Update the meshes $\T_{h,\ell}$ and $\T^2_{h,\ell}$ and recompute the error indicators.
    \end{enumerate}
\end{enumerate}

	\subsection{Numerical results}
	
	In this section, we present a sample test for the elliptic interface problem in the FD-DLM formulation, focusing on both a priori and a posteriori error analysis. We also discuss some numerical results regarding the preconditioners. 
	
\begin{figure}
		\centering
		\subfloat[\label{fig:circle_geo}]{\includegraphics[width=0.2\linewidth]{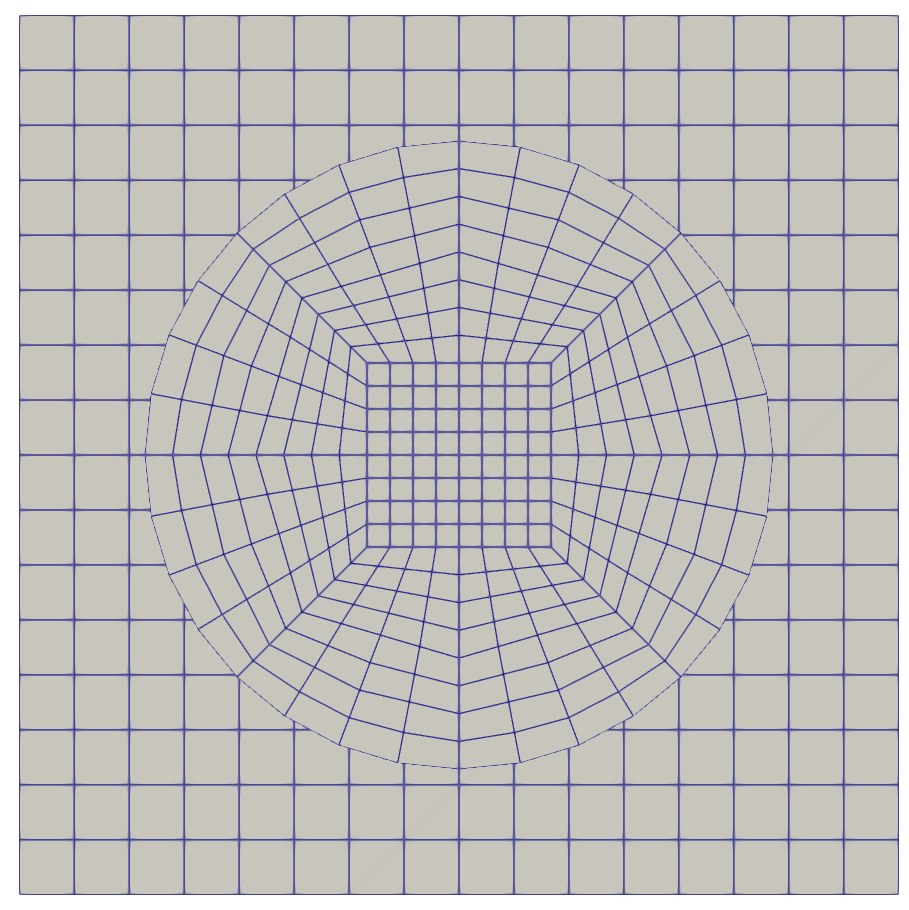}}\qquad\qquad
		\subfloat[\label{fig:circle_geo2}]{\includegraphics[width=0.2\linewidth]{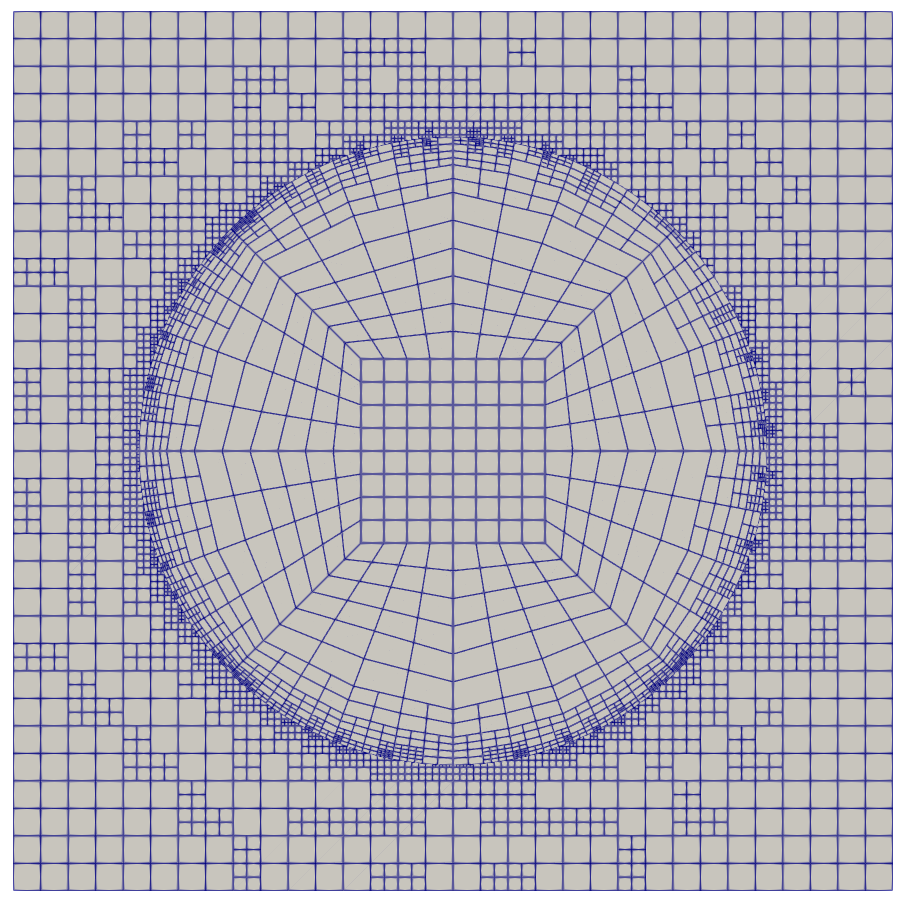}}\qquad\qquad
		\subfloat[\label{fig:circle_sol}]{\includegraphics[width=0.25\linewidth]{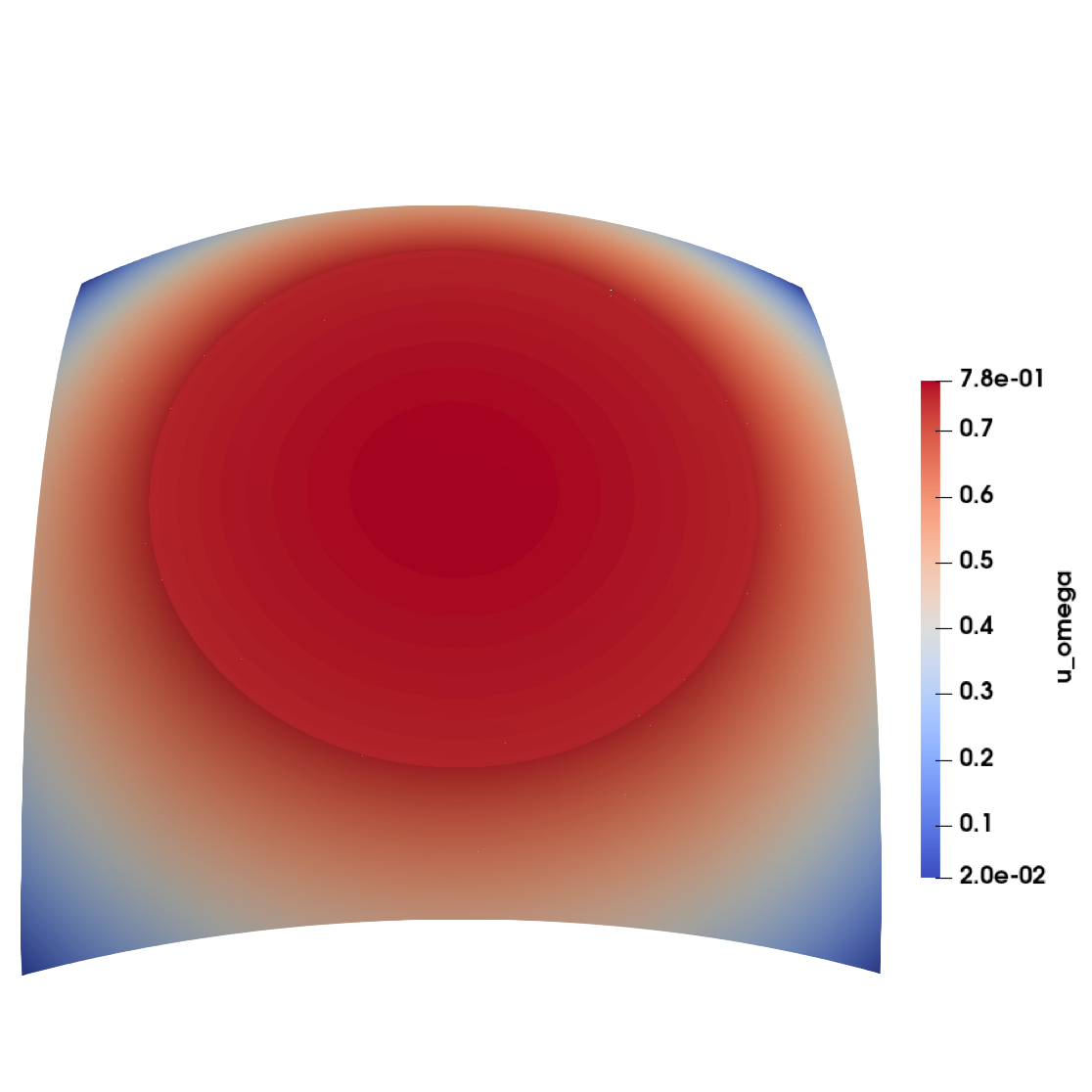}}\quad
		\caption{\textsc{(a)} Domain configuration of the immersed circle, initial mesh. \textsc{(b)}  Domain configuration after adaptive refinement. \textsc{(c)}  Profile of the computed solution.}
		\label{fig:examples}
	\end{figure}

	We consider the example of the circle $\Omega_2 =\{\x\in\R^2:\,|\x|\le1\}$ immersed in the square $\Omega=[-1.4,1.4]^2$ already discussed in~\cite{najwa-posteriori}. Other tests cases can be found in~\cite{najwa,najwa-posteriori}. In this case, the following exact solution is available, and for $(\nu, \nu_2)=(1,10)$ its analytic expression is
	\begin{equation*}
		u_1(x,y) = \dfrac{4-x^2-y^2}{4},
		\qquad
		u_2(x,y)=\dfrac{31-x^2-y^2}{40}.
	\end{equation*}
	The geometric configuration, the initial mesh and the last adaptive refinement together with the profile of the extended solution $u_h$ are reported in Figure~\ref{fig:examples}. 
	
	We discretize the problem by the $(\Qoly_1,\Qoly_1+\Bbub,\Poly_0)$ element. The coupling term is the scalar product in $\LdOdd$ and its construction is carried out by exact integration. 
	The problem is solved with both uniform mesh refinement and the adaptive strategy presented in Section~\ref{sec:adaptive}, guided by error indicators defined in~\eqref{eq:eta}. The refinement fraction is $\alpha=0.6$. As shown in \cite{najwa-posteriori}, this refinement fraction consistently yields optimal convergence rates.

\begin{figure}
		\includegraphics[width=0.38\linewidth]{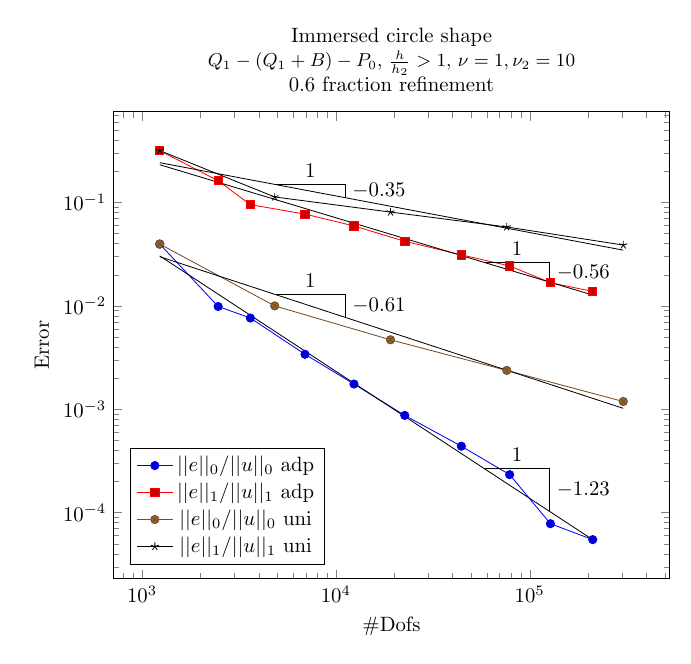}\qquad
		\includegraphics[width=0.38\linewidth]{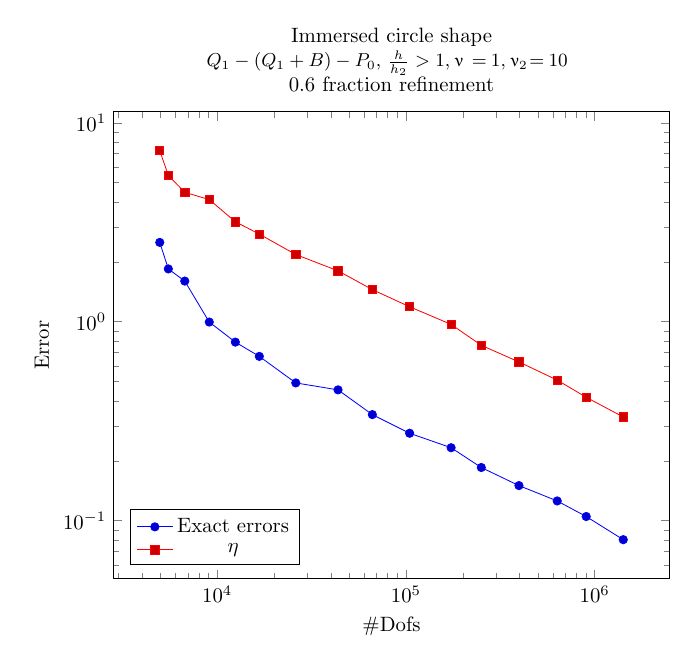} 
		\caption{Left: comparison of the error $e=u-u_h$ between uniform (uni) and adaptive (adp) refinement. Right: comparison between the exact error $\|e\|_{1}$ and the a posteriori indicator $\eta$.}
		\label{fig:rate_mesh_ratio}
	\end{figure}

	We examine the convergence rates of the error measured in the $L^2$ and $H^1$ norms of the solution $u$. The convergence behaviors are illustrated in Figure~\ref{fig:rate_mesh_ratio} (left), comparing errors against the number of degrees of freedom for uniform and adaptive mesh refinements. Uniform refinement confirms convergence rates of $O(h)$ for the $L^2$ norm and $O(h^{1/2})$ for the $H^1$ norm, consistent with theoretical results based on the regularity of the solution. In contrast, the adaptive refinement process significantly enhances the convergence rates, achieving optimal rates of $O(1/\sharp\mathrm{Dofs})$ for the $L^2$ norm and $O(1/\sharp\mathrm{Dofs}^{1/2})$ for the $H^1$ norm. This highlights the effectiveness of adaptive refinement in capturing localized solution irregularities.
	To assess the reliability of the error estimator, we compare the computed error estimator against the actual error. The results, reported in Figure~\ref{fig:rate_mesh_ratio} (right), reveal a consistent linear relationship between the error estimator and the exact error. The ratio between them remains constant, confirming reliability and efficiency of the estimator in predicting the error across all mesh refinement levels.
	The behavior of our scheme is independent of the magnitude of the coefficient jump and of mesh sizes ratio. The interested reader can refer to~\cite{najwa-posteriori}.
	
	\subsubsection{Preconditioning strategies} 
	In this section, we choose both definitions of the coupling space, namely $\Lambda = H^1(\Omega_2)$ and $\Lambda = \Vtd$. We consider again the example of the circular interface with coefficients $(\nu, \nu_2) = (1, 10)$, as done in~\cite{alshehri2025multigrid}. For other test cases, we refer to~\cite{boffi2022parallel}, where we used the preconditioners for solving fluid-structure interaction problems.
	
	We analyze two discretization schemes: the first one is based on the $(\Qoly_1, \Qoly_1, \Qoly_1)$ element and employs both the $H^1$ and $L^2$ coupling terms, the latter is given by the $(\Qoly_1, \Qoly_1 + \Bbub, \Poly_0)$ element and the coupling term is the $L^2$ scalar product.
	
	The linear system is solved using GMRES with a tolerance of $10^{-12}$. Moreover, we employ the preconditioners defined in~\eqref{eq:prec}. Their action requires the inversion of the two diagonal blocks. This can be done by using a direct solver, denoted by ``$\mathsf{d}$'', or by a multigrid one, ``$\mathsf{m}$''. We consider three configurations for inverting the blocks of the preconditioner matrix: $\dd$, $\md$, and $\mm$. Here, $\md$ means that the first block is inverted with multigrid, while $\dd$ and $\mm$ refer to the inversion of both blocks with either the direct or multigrid method.
	
	Table~\ref{tbl:H1cont_iterations} reports the iteration count $its$ and solution time $T_{sol}$ for the case $\Lambda = H^1(\Omega_2)$ with the element $(\Qoly_1, \Qoly_1, \Qoly_1)$. The results indicate that the diagonal preconditioner combined with $\dd$ inversion performs poorly. However, the $\md$ and $\mm$ configurations are both effective in controlling the iteration count and solution time, with $\md$ showing a slightly better performance. For the lower triangular preconditioner, all three inversion strategies yield comparable performance.

	Table~\ref{tbl:L2cont_iterations} presents the results for $\Lambda = \Vtd$ using the same element. In this case, the diagonal preconditioner performs best when inverted using the $\dd$ strategy with solution time comparable among all cases. On the other hand, all three strategies demonstrate similar overall performance for the triangular preconditioner.

	Lastly, Table~\ref{tbl:L2disc_iterations} presents results for $\Lambda = \Vtd$ with the element $(\Qoly_1, \Qoly_1 + \Bbub, \Poly_0)$. In this setting, $\mm$ shows poor performance. In contrast, the $\dd$ and $\md$ strategies yield comparable results in terms of both iteration count and solution time when applied with the triangular preconditioner.

	\begin{table}
		\begin{center}
			\resizebox{\textwidth}{!}{
			\begin{tabular}{r|r|r|r|r|r|r|r|r|r|r|r|r}
				\hline
				\multicolumn{13}{c}{{\bf Circular immersed domain -- Mesh refinement test}} \\
				\multicolumn{13}{c}{{\textit{Exact integration of the coupling term, $\Lambda=H^1(\Omega_2)$, $(\Qoly_1,\Qoly_1,\Qoly_1)$}}} \\
				\hline
				$dofs$    & \multicolumn{6}{c|}{$\mathsf{P}_{\mathsf{diag}}$} 		& \multicolumn{6}{c}{$\mathsf{P}_{\mathsf{tri}}$} \\
				& \multicolumn{2}{c|}{$\dd$}   & \multicolumn{2}{c|}{$\md$} 		& \multicolumn{2}{c|}{$\mm$} 
				& \multicolumn{2}{c|}{$\dd$}   & \multicolumn{2}{c|}{$\md$} 		& \multicolumn{2}{c}{$\mm$}\\
				& $its$           & $T_{sol}(s)$ & $its$           & $T_{sol}(s)$	& $its$           & $T_{sol}(s)$           & $its$           & $T_{sol}(s)$& $its$           & $T_{sol}(s)$& $its$           & $T_{sol}(s)$ \\
				\hline
				755     & 16 & 9.49e--03& 68 & 9.00e--03& 77 & 3.00e--01 & 10 & 4.34e--03& 33 & 3.09e--02& 34 & 5.47e--03\\
				2,915    &22 & 2.60e--02& 82 & 4.86e--02& 90 & 1.40e+00   & 23 & 2.49e--02& 42 & 6.31e--02& 43 & 2.75e--02\\
				11,459   &539 & 6.83e+07& 84 & 4.21e--01& 90 & 5.46e+00  & 38 & 3.41e--01& 44 & 5.85e--01& 45 & 2.12e--01\\
				45,443   &-&-&86&2.32e+00& 99 & 2.46e+01                & 46 & 2.48e+00& 47 & 4.13e+00& 47 & 1.25e+00\\
				180,995  &-&-&86&1.15e+01& 99 & 1.01e+02                & 43 & 1.20e+01& 46 & 3.52e+01& 46 & 6.32e+00\\
				722,435  &-&-&86&6.34e+01& 101 & 4.21e+02               & 50 & 7.75e+01& 51 & 7.19e+02& 49 & 3.67e+01\\
				\hline
			\end{tabular}}
			\caption{Mesh refinement test for the immersed circular shape, $H^1$ coupling term with mesh intersection and $(\Qoly_1,\Qoly_1,\Qoly_1)$ element. The simulations are run on one processor. Legend: $dofs$ = degrees of freedom; $its$ = GMRES iterations; $T_{sol}$ = CPU time to solve the linear system. All CPU times are reported in seconds. }
			\label{tbl:H1cont_iterations}
		\end{center}
	\end{table}
	\begin{table}
		\begin{center}
			\resizebox{\textwidth}{!}{
			\begin{tabular}{r|r|r|r|r|r|r|r|r|r|r|r|r}
				\hline
				\multicolumn{13}{c}{{\bf Circular immersed domain -- Mesh refinement test}} \\
				\multicolumn{13}{c}{{\textit{Exact integration of the coupling term, $\Lambda=\Vtd$, $(\Qoly_1,\Qoly_1,\Qoly_1)$}}} \\
				\hline
				$dofs$    & \multicolumn{6}{c|}{$\mathsf{P}_{\mathsf{diag}}$} 		& \multicolumn{6}{c}{$\mathsf{P}_{\mathsf{tri}}$} \\
				& \multicolumn{2}{c|}{$\dd$}   & \multicolumn{2}{c|}{$\md$} 		& \multicolumn{2}{c|}{$\mm$} 
				& \multicolumn{2}{c|}{$\dd$}   & \multicolumn{2}{c|}{$\md$} 		& \multicolumn{2}{c}{$\mm$}\\
				& $its$           & $T_{sol}(s)$ & $its$           & $T_{sol}(s)$	& $its$           & $T_{sol}(s)$           & $its$           & $T_{sol}(s)$& $its$           & $T_{sol}(s)$& $its$           & $T_{sol}(s)$ \\
				\hline
				755     & 19 & 3.02e--03& 84 & 1.14e--02& 84 & 3.21e--01		& 11 & 2.72e--03& 40 & 1.25e--02& 49 & 3.80e--01 \\
				2,915    & 48 & 3.73e--02& 103 & 6.69e--02& 104 & 1.57e+00		& 21 & 2.52e--02 & 49 & 5.55e--02& 54 & 1.67e+00 \\
				11,459   & 82 & 4.24e--01& 110 & 4.61e--01& 112 & 6.68e+00	    & 43 & 3.37e--01& 54 & 3.99e--01& 59 & 7.13e+00\\
				45,443   & 68 & 2.22e+00 & 118 & 3.14e+00& 120 & 2.97e+00	& 55 & 2.92e+00& 59 & 3.04e+00& 74 & 3.76e+01\\
				180,995  & 46 & 7.66e+00& 147 & 1.91e+01& 149 & 1.52e+02		& 60 & 1.56e+01& 78 & 2.01e+01& 78 & 1.57e+02 \\
				722,435  & 38 & 3.50e+01&-&-&-&-								& 79 & 4.82e+02& 79 & 1.07e+02& 81 & 6.64e+02\\
				\hline
			\end{tabular}}
			\caption{Mesh refinement test for the immersed circular shape, $L^2$ coupling term with mesh intersection and $(\Qoly_1,\Qoly_1,\Qoly_1)$ element. The simulations are run on one processor. Same format as in Table \ref{tbl:H1cont_iterations}.}
			\label{tbl:L2cont_iterations}
		\end{center}
	\end{table}
	\begin{table}
		\begin{center}
			\resizebox{\textwidth}{!}{
			\begin{tabular}{r|r|r|r|r|r|r|r|r|r|r|r|r}
				\hline
				\multicolumn{13}{c}{{\bf Circular immersed domain -- Mesh refinement test}} \\
				\multicolumn{13}{c}{{\textit{Exact integration of the coupling term, $\Lambda=\Vtd$, $(\Qoly_1,\Qoly_1+\Bbub,\Poly_0)$}}} \\
				\hline
				$dofs$    & \multicolumn{6}{c|}{$\mathsf{P}_{\mathsf{diag}}$} 		& \multicolumn{6}{c}{$\mathsf{P}_{\mathsf{tri}}$} \\
				& \multicolumn{2}{c|}{$\dd$}   & \multicolumn{2}{c|}{$\md$} 		& \multicolumn{2}{c|}{$\mm$} 
				& \multicolumn{2}{c|}{$\dd$}   & \multicolumn{2}{c|}{$\md$} 		& \multicolumn{2}{c}{$\mm$}\\
				& $its$           & $T_{sol}(s)$ & $its$           & $T_{sol}(s)$	& $its$           & $T_{sol}(s)$           & $its$           & $T_{sol}(s)$& $its$           & $T_{sol}(s)$& $its$           & $T_{sol}(s)$ \\
				\hline
				1,058 	& 19 & 2.81e--03& 80 & 9.90e--03& 107 & 2.53e--01 		& 9 & 2.04e--03& 39 & 7.84e--03& 60 & 2.85e--01 \\
				4,162 	& 52 & 3.14e--02& 103 & 5.82e--02& 137 & 1.30e+00		& 20 & 1.74e--02& 49 & 4.30e-02& 77 & 1.45e+00 \\
				16,514 	& 82 & 2.48e--01& 111 & 3.03e--01& 325 & 1.17e+01		& 41 & 1.84e--01& 54 & 2.43e--01& 111 & 8.17e+00\\
				65,794 	& 80 & 1.43e+00& 118 & 1.97e+00&-&-					& 59 & 1.71e+00& 73 & 2.20e+00& 989 & 2.90e+02\\
				262,658 	& 76 & 7.35e+00& 145 & 1.35e+01&-&-					& 77 & 1.16e+01& 77 & 1.16e+01&-&- \\
				1,049,602 & 52 & 2.64e+01&-&-&-&-								& 81 & 6.30e+01& 81 & 6.33e+01&-&-\\
				\hline
			\end{tabular}}
			\caption{Mesh refinement test for the immersed circular shape, $L^2$ coupling term with mesh intersection and $(\Qoly_1,\Qoly_1+\Bbub,\Poly_0)$ element. The simulations are run on one processor. Same format as in Table \ref{tbl:H1cont_iterations}.}
			\label{tbl:L2disc_iterations}
		\end{center}
	\end{table}

	\section{Application to fluid-structure interaction problems}\label{sec:fsi}
	
	In this section we apply our fictitious domain formulation for modeling fluid-structure interaction problems where a viscous hyper-elastic solid body is immersed in a Newtonian fluid.
	
	\subsection{The continuous problem}
	
	We study fluid-structure interaction problems in a fixed domain $\Omega\subset\R^d$ with $d=2,3$. Such domain is decomposed into two time dependent regions, $\Oft$ and $\Ost$, occupied by the evolving fluid and solid, respectively, at time instant $t$. We assume that $\Oft$ and $\Ost$ have codimension zero, even if more general situations, such as immersed thin structures, can be modeled by our approach (see e.g.~\cite{2015,annese,ALZETTA,heltaizunino} for more details). Moreover, for simplicity of exposition, we assume that the immersed interface $\Gamma_t = \overline{\Omega}_t^s\cap\overline{\Omega}_t^f$ cannot touch $\partial\Omega$, that is $\Gamma_t\cap\partial\Omega=\emptyset$.
	
	We describe the fluid dynamics in Eulerian framework and we denote by $\x$ the associated variable. The dynamics of incompressible Newtonian fluids is governed by the Navier--Stokes equations in~$\Oft$. Given the viscosity $\nu_f>0$, the density $\rho_f$, the velocity $\u_f$ and the pressure $p_f$, the Cauchy stress tensor for viscous fluids reads
	\begin{equation}
		\ssigma_f = -p_f\I + \nu_f\Grads(\u_f),
	\end{equation}
	where the symbol $\Grads(\cdot)$ refers to the symmetric gradient operator, that is $\Grads(\cdot) = (\Grad\cdot+\Grad^\top\cdot)/2$. 
	
\begin{figure}
		\centering
		\includegraphics[width=0.6\linewidth]{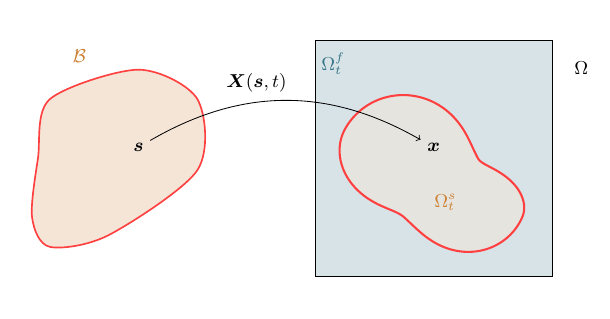}
		\caption{Geometric configuration of the FSI problem. The domain $\Omega$ is fixed in time and behaves as a container. The fluid and solid regions, $\Oft$ and $\Ost$ respectively, are time dependent. The Lagrangian description is employed to represent the solid deformation: the reference domain $\B$ is mapped into the actual position $\Ost$ through the map $\X$.}
		\label{fig:geo}
	\end{figure}

	The solid deformation is described in Lagrangian setting. We consider a fixed domain $\B\in\R^d$ playing the role of reference domain for the solid configuration. At each time instant $t$, the actual configuration of the solid is then obtained through the action of the deformation map ${\X(\cdot,t):\B\rightarrow\Ost}$. We have that
	\begin{equation}
		\x = \X(\s,t)\qquad\text{for }\x\in\Ost,
	\end{equation}
	where the Lagrangian variable is denoted by $\s\in\B$. An example of geometric configuration for the FSI system is sketched in Figure~\ref{fig:geo}.
	
	The immersed solid body is made of incompressible viscous hyper-elastic material. Thus, the associated Cauchy stress tensor $\ssigma_s$ takes into account both the viscous and hyper-elastic constitutive laws. Indeed, we have $\ssigma_s = \ssigma_s^v + \ssigma_s^e$. The viscous contribution $\ssigma_s^v$ is written similarly to $\ssigma_f$. Denoting by $\nu_s>0$ and~$\rho_s$ the solid viscosity and density, respectively, and given the solid velocity $\u_s$ we have
	\begin{equation}
		\ssigma_s^v = -p_s\I + \nu_s\Grads(\u_s),
	\end{equation}
	where the pressure $p_s$ plays the role of Lagrange multiplier to enforce the incompressibility of the solid material. Concerning the elastic contribution, $\ssigma_s^e$ is expressed in terms of the first Piola--Kirchhoff elasticity tensor $\PP$, which gives
	\begin{equation}
		\ssigma_s^e = |\F|^{-1}\PP\F.	
	\end{equation}
	More precisely, $\F=\grads\X$ denotes the deformation gradient, and $|\F|$ its determinant. The incompressibility condition implies that $\F$ is constant in time. Moreover, $|\F|=1$ if the initial configuration $\Os_0$ coincides with $\B$.
	
	The motion of the immersed structure is represented by the following kinematic condition, which relates the solid velocity $\u_s$ to the time derivative $\X$,
	\begin{equation}\label{eq:motion}
		\u_s(\x) = \derivative{\X(\s,t)}{t}\qquad\text{for }\x\in\Ost.
	\end{equation}

	Moreover, fluid and solid models are coupled by means of transmission conditions which enforce continuity of velocities and stresses along the interface $\Gamma_t$.

	The mathematical model governing our system is the following
	\begin{equation*}
		\begin{aligned}
			& \rho_f \left(\derivative{\u}{t}+\u\cdot\Grad\u\right) = \div\ssigma_f&&\text{in }\Oft\\
			& \div\u_f = 0&&\text{in }\Oft\\
			& \rho_s\frac{\partial^2\X}{\partial t^2} = \div_\s\left(|\F|\ssigma_s^v\F^{-\top}+\PP(\F)\right)&&\text{in }\B\\
			& \div\u_s = 0&&\text{in }\Ost\\
			& \u_f = \derivative{\X}{t}&&\text{on }\Gamma_t\\
			& \ssigma_f\n_f = -\left(\ssigma_s^v+|\F|^{-1}\PP\F^{\top}\right)\n_s&&\text{on }\Gamma_t.
		\end{aligned}
	\end{equation*}
	The model is then completed by adding the no-slip condition $\u_f=0$ on $\partial\Omega$ and the initial conditions
	\begin{equation*}
		\u_f(0) = \u_{f,0} \quad\text{ in }\Of_0,
		\qquad\qquad
		\u_s(0) = \u_{s,0} \quad\text{ in }\Os_0,
		\qquad\qquad
		\X(0) = \X_0  \quad\text{ in }\B.
	\end{equation*}

	The first step towards the derivation of the fictitious domain formulation with distributed Lagrange multiplier consists in defining velocity and pressure on the entire $\Omega$, by extending them to the solid region. We define the extended variables as
	\begin{equation*}
		\u = \begin{cases}
			\u_f & \text{in }\Oft\\
			\u_s & \text{in }\Ost,
		\end{cases}
		\qquad
		p = \begin{cases}
			p_f & \text{in }\Oft\\
			p_s & \text{in }\Ost.
		\end{cases}
	\end{equation*}
	In order to have equivalence with the original problem, the kinematic condition~\eqref{eq:motion} becomes a constraint. Indeed, we impose that the velocity $\u$ must equal $\partial\X/\partial t$ on the fictitiously extended region, that is
	\begin{equation}\label{eq:constraint}
		\u(\X(\s,t),t) = \derivative{\X(\s,t)}{t}\qquad\text{for }\x=\X(\s,t).
	\end{equation}
	This constraint is enforced at variational level by means of a distributed Lagrange multiplier and a coupling bilinear form $\c(\cdot,\cdot)$, both defined on a suitable functional space $\LL$, as explained in the case of elliptic interface problems (see~\eqref{eq:coupling_def}). Also in this case, $\LL$ is defined either as $\Hubd$, namely the dual space of $\Hub$, or $\Hub$ itself. In the former case, $\c(\cdot,\cdot)$ is the duality pairing between $\Hub$ and $\Hubd$, whereas in the latter $\c(\cdot,\cdot)$ is defined as the scalar product in $\Hub$.
	
	The fluid-structure interaction problem in fictitious domain framework is the following.
	
	\begin{pro}
		\label{pro:problem_fictitious}
		Given $\u_0\in\Huo$ and $\X_0:\B\rightarrow\Omega$, $\forall t\in (0,T)$, find $\u(t)\in\Huo$, $p(t)\in\Ldo$, $\X(t)\in \mathbf{W}^{1,\infty}(\B)$ and $\llambda(t)\in \LL$, such that
		\begin{equation*}
			\begin{aligned}
				&\notag \rho_f\left(\derivative{\u(t)}{t},\v\right)_\Omega+b(\u(t),\u(t),\v)+a(\u(t),\v)\\
				&\hspace{2.5cm}-(\div\v,p(t))_\Omega+\c(\llambda(t),\v(\X(t)))=0 
				&&\quad\forall\v\in\Huo\\
				&(\div\u(t),q)_\Omega=0  &&\quad\forall q\in \Ldo\\
				& \label{eq:fict_solid} \dr \bigg(\frac{\partial^2\X}{\partial t^2},\Y \bigg)_\B+(\PP(\F),\grads\Y )_\B-\c(\llambda(t),\Y)=0&&\quad\forall\Y\in \Hub \\
				& \c\bigg(\mmu,\u(\X(\cdot,t),t)-\derivative{\X(\cdot,t)}{t} \bigg)=0 &&\quad\forall\mmu\in\LL \\
				&\u(\x,0)=\u_0(\x) &&\quad\forall\x\in\Omega\\
				&\X(\s,0)=\X_0(\s) &&\quad\forall\s\in\B.
			\end{aligned}
		\end{equation*}
	\end{pro}

	The reader interested in the detailed derivation of the above problem may refer to~\cite{2015}. We define the following bilinear and trilinear forms on the fluid domain
	\begin{equation*}
		a(\u,\v)=\left(\nu\Grads(\u),\Grads(\v)\right)_\Omega,
		\qquad
		b(\u,\v,\w)=\frac{\rho_f}{2}\left((\u\cdot\nabla\v,\w)_\Omega-(\u\cdot\nabla\w,\v)_\Omega\right),
	\end{equation*}
	where the extended viscosity $\nu$ equals $\nu_f$ in $\Oft$ and $\nu_s$ in $\Ost$. Moreover, $\dr = \rho_s-\rho_f$.

	\begin{rem}
		Existence and uniqueness of the solution for Problem~\ref{pro:problem_fictitious} have been proved in~\cite{2020} for a simplified version of the problem by following the Galerkin approximation technique used in~\cite{temam}. In~\cite{2020}, a linearized version of the problem is considered: indeed the convective term is neglected and a linear constitutive law is chosen for the solid, that is $\PP(\F) = \kappa\F$. Moreover, the coupling form $\c(\cdot,\cdot)$ is set to the scalar product in $\Hub$. The following space-time regularity is ensured
		\begin{gather*}
			\u \in \mathbf{L}^\infty(0,T;\mathscr{V}_0) \cap \mathbf{L}^2(0,T;\mathscr{H}_0),
			\qquad
			p \in L^2(0,T;\Ldo),
			\qquad
			\llambda \in \mathbf{L}^2(0,T;\Hub),\\
			\X \in \mathbf{L}^\infty(0,T;\Hub) \text{ with } \derivative{\X}{t}\in \mathbf{L}^\infty(0,T;\LdBd) \cap \mathbf{L}^2(0,T;\Hub),
		\end{gather*}
		where
		\begin{gather*}
		\mathscr{V} = \{ \v\in(\mathcal{D}(\Omega))^d:\div\v=0\},\\
		\mathscr{H}_0 = \text{ the closure of }\mathscr{V}\text{ in }\Huo,
		\qquad
		\mathscr{V}_0 = \text{ the closure of }\mathscr{V}\text{ in }\mathbf{L}_0^2(\Omega),
		\end{gather*}
		and $\mathcal{D}(\Omega)$ denotes the space of infinitely differentiable functions with compact support in $\Omega$.
	\end{rem}

	\subsection{Discretization}
	
	We discretize Problem~\ref{pro:problem_fictitious} by mixed finite elements.
	
	Similarly to what we did for the elliptic interface problem, fluid and solid (reference) domain are treated independently. We thus partition $\Omega$ by a mesh $\T_h^\Omega$ with size $h_\Omega$, while we decompose the domain $\B$ by a mesh $\T_h^\B$ with size $h_\B$. We emphasize that the considered meshes are fix throughout the entire evolution of the system.
	
	We then choose four finite dimensional subspaces, namely ${\Vline_h\subset\Huo}$, ${Q_h\subset\Ldo}$, ${\Sline_h\subset\Hub}$ and ${\LL_h\subset\LL}$. In order to design a stable method, the velocity and pressure spaces $(\Vline_h,Q_h)$ must be a compatible pair for the Stokes problem~\cite{mixedFEM}. Compatibility requirements involving the solid discrete spaces $\Sline_h$ and $\LL_h$ have been previously established for the elliptic interface problems and apply straightforwardly to this case. For more details, see Section~\ref{sec:discrete_elliptic} and discrete spaces therein.
	
	As time advancing scheme, we choose Backward Euler on a uniform partition of the time domain $[0,T]$. We denote the time step by $\dt$ and the time grid by~$\{t_n\}_{n=0}^N$ with $t_n=n\dt$. For a generic function $v$, we write $v^n=v(t_n)$ and
	\begin{equation*}
		\derivative{v}{t}(t_{n+1})\approx\frac{\v^{n+1}-\v^{n}}{\dt},
		\qquad
		\frac{\partial^2 v}{\partial t^2}(t_{n+1})\approx\frac{\v^{n+1}-2\v^{n}+\v^{n-1}}{\dt^2}\text{.}
	\end{equation*}

	The fully discrete problem reads as follows.
	
	\begin{pro}
		\label{pro:discrete}
		Given $\u_h^0\in\Vline_h$ and $\X_h^0\in\Sline_h$, for $n=1,\dots,N$ find $u_h^n\in\Vline_h$, $p_h^n\in Q_h$, $\X^n_h\in\Sline_h$, and $\llambda^n_h\in\LL_h$, such that
		\begin{equation*}
			\begin{aligned}
				&\notag\rho_f\bigg(\frac{\u^{n+1}_h-\u^{n}_h}{\dt},\v\bigg)+b(\u^{n}_h,\u^{n+1}_h,\v_h)+a(\u^{n+1}_h,\v_h)\\
				&\hspace*{3.75cm}-(\div\v_h,p^{n+1}_h)+\c(\llambda^{n+1}_h,\v_h(\X^{n}_h))=0&&\forall\v_h\in\Vline_h\\
				&(\div\u^{n+1}_h,q_h)=0  &&\forall q_h\in Q_h\\
				&\notag \dr \bigg(\frac{\X^{n+1}_h-2\X^n_h+\X^{n-1}_h}{\dt^2},\Y_h \bigg)_\B+(\PP(\F^{n+1}),\grads\Y_h )_\B\\
				&\hspace*{7cm}-\c(\llambda^{n+1}_h,\Y_h)=0 &&\forall\Y_h\in \Sline_h\\
				& \label{eq:initiTime}\c\bigg(\mmu_h,\u^{n+1}_h(\X^{n}_h)-\frac{\X_h^{n+1}-\X_h^{n}}{\dt} \bigg)=0  &&\forall\mmu_h\in\LL_h.
			\end{aligned}
		\end{equation*}
	\end{pro}

	We observe that at the first time step, we need the value $\X_h^{-1}$ to initialize the advancing scheme. This value can be obtained by solving the following equation involving the initial conditions $\u_{s,0}$ and $\X_h^0$
	\begin{equation*}
		\u_{s,0} = \frac{\X_h^0 - \X_h^{-1}}{\dt} \qquad \text{in }\B.
	\end{equation*}
	
	The solution of a fully implicit scheme would require massive computational resources for dealing with several nonlinearities; for this reason, some terms are treated explicitly. This is the case of the nonlinear convective term $b(\cdot,\cdot,\cdot)$, which involves the velocity $\u_h^n$ at the previous time instant. The same idea applies to the coupling terms $\c(\llambda^{n+1}_h,\v_h(\X^{n}_h))$ and $\c(\mmu_h,\u_h(\X^{n}_h))$, where the actual position of the immersed structure is taken into account by composing~$\u_h$ and~$\v_h$ with~$\X^{n}_h$ instead of~$\X^{n+1}_h$. Moreover, the term $\PP(\F^{n+1})$ may require additional linearization according to the involved constitutive law.
	
	In any case, our semi-implicit scheme is unconditionally stable, as stated by the following proposition.
	
	\begin{prop}[\cite{2015}]
		The discrete Problem~\ref{pro:discrete} is unconditionally stable in time, without any restriction on the choice of time step $\dt$.
	\end{prop}

	\begin{rem}
		We chose Backward Euler for the time discretization, but high order schemes may also be considered in the spirit of~\cite{dong10,okamoto10,chen13}. In~\cite{wolf}, we proved that the second order backward differentiation formula \textsf{BDF2} gives an unconditionally stable scheme. The same result is also true for the Crank--Nicolson scheme implemented with the midpoint rule. A proof of unconditional stability is not available for the trapezoidal Crank--Nicolson, but several numerical tests did not show any kind of instability.
	\end{rem}

	\fc
	\begin{rem}
		The incompressibility of both fluid and solid is handled by imposing the divergence-free constraint on the extended velocity $\u$, and guarantees conservation of mass at the continuous level. Conversely, at the discrete level, the divergence-free constraint may not be exactly imposed. Indeed, it is well-known that many popular mixed finite element pairs for the Stokes equation (Hood--Taylor, Bercovier--Pironneau, MINI, etc., see~\cite{mixedFEM}) conserve mass only in approximate way~\cite{linke}, while the construction of conforming divergence-free finite elements is not trivial and requires special recipes (see~\cite{arnold2006finite,cockburn2007note,da2023virtual}, for instance). The mass conservation properties of our approach in terms of conforming finite elements have been discussed in~\cite{boffi2013mass,CABG25} and further details will be addressed in forthcoming works.
	\end{rem}
	\cf
	
	We rewrite Problem~\ref{pro:discrete} in matrix form by considering the special case of a linear constitutive law for the solid, that is $\PP(\F)=\kappa\F$. In the following, $\boldsymbol{\phi}$, $\psi$, $\boldsymbol{\chi}$, $\boldsymbol{\zeta}$ denotes the basis functions for $\Vline_h$, $Q_h$, $\Sline_h$ and $\LL_h$, respectively. We have
	
	\begin{equation}\renewcommand{\arraystretch}{1.5}
		\label{eq:matrix_fsi}
		\left[\begin{array}{@{}ccc|c@{}}
			\Amatr_f		& -\Bmatr^\top	& 0			& \Cmatr_f(\X_h^n)^\top \\
			-\Bmatr		& 0	& 0 			& 0 \\
			0 		& 0 	& \Amatr_s   		& -\Cmatr_s^\top \\
			\hline
			\Cmatr_f(\X_h^n) 	& 0   	& -\frac{1}{\dt}\Cmatr_s 	& 0 \\
		\end{array}\right]
		\left[ \begin{array}{c}
			\u_h^{n+1} \\
			p_h^{n+1} \\
			\X_h^{n+1} \\
			\hline
			\llambda_h^{n+1} \\
		\end{array}\right]
		= 
		\left[ \begin{array}{c}
			\mathsf{g}_1 \\
			0 \\
			\mathsf{g}_2 \\
			\hline
			\mathsf{g}_3 
		\end{array}\right],
	\end{equation}\renewcommand{\arraystretch}{1.15}
	where
	\begin{equation*}
		\begin{aligned}
			& \Amatr_f = \frac{\rho_f}{\dt} \Mmatr_f + \Kmatr_f,
			\qquad
			(\Mmatr_f)_{ij} = \left( \boldsymbol{\phi}_j,
			 \boldsymbol{\phi}_i \right)_\Omega,
			\qquad
			(\Kmatr_f)_{ij} = a\left(\boldsymbol{\phi}_j, \boldsymbol{\phi}_i\right)+b(\u_h^n,\boldsymbol{\phi}_j,
			\boldsymbol{\phi}_i), \\
			& \Bmatr_{ki} = \left( \div \boldsymbol{\phi}_i, \psi_k\right)_\Omega,\\
			& (\Amatr_s)_{ij} = \frac{\dr}{\dt^2} \Mmatr_s + \Kmatr_s,
			\qquad
			(\Mmatr_s)_{ij} =  \left( \boldsymbol{\chi}_j, \boldsymbol{\chi}_i \right)_\B,
			\qquad
			(\Kmatr_s)_{ij} = \kappa \left( \nabla_s \boldsymbol{\chi}_j, \nabla_s \boldsymbol{\chi}_i \right)_\B,\\
			& (\Cmatr_f(\X_h^n))_{\ell j} = \c(\boldsymbol{\zeta}_\ell, \boldsymbol{\phi}_j(\X_h^n)),\\
			& (\Cmatr_s)_{\ell j} = \c(\boldsymbol{\zeta}_\ell, \boldsymbol{\chi}_j),\\
			& \mathsf{g}_1 = \frac{\rho_f}{\dt}\Mmatr_f\u_h^n,
			\quad \mathsf{g}_2 = \frac{\dr}{\dt^2} \Mmatr_s(2\X_h^n-\X_h^{n-1}),
			\quad \mathsf{g}_3 = -\frac{1}{\dt} \Cmatr_s \X_h^n.
		\end{aligned}
	\end{equation*}

	\begin{figure}
		\centering
		\includegraphics[width=0.5\linewidth]{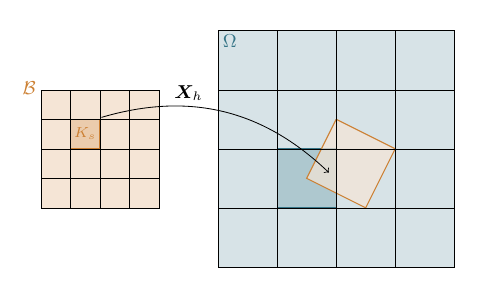}
		\caption{Superimposition of a solid element $\E_s\in\T_h^\B$ on the background mesh~$\T_h^\Omega$ through the action of the Lagrangian map~$\X$.}
		\label{fig:mapping}
	\end{figure}

	The coupling matrix $\Cmatr_f(\X_h^n)$ corresponds to the terms of kind $\c(\mmu_h,\v_h(\X_h^n))$ and represents the interaction between the fluid and the immersed solid. We point out that its construction requires the integration over $\B$ of $\v_h$ (defined on $\T_h^\Omega$) composed with the deformation map $\X_h^n$: this takes care of the actual position of the solid body with respect to the background mesh, as sketched in Figure~\ref{fig:mapping}. It is clear that the assembly of $\Cmatr_f(\X_h^n)$ falls in the framework we described in Section~\ref{sec:coupling}.
	
	In the case of a nonlinear solid constitutive law, the matrix $\Amatr_s$ depends on the deformation $\X$, thus the system must be solved by employing a solver for nonlinear systems of equations, such a fixed point iteration or a Newton-like method.

	At each time step, system~\eqref{eq:matrix_fsi} can be interpreted as a stationary saddle point problem, which admits a unique solution. Stability and well-posedness of the finite element discretization have been proved in~\cite{stat}: the method gives optimal convergence rates, which depend on the regularity of the solution. Well-posedness and convergence properties in case of inexact integration of the coupling term have been discussed in~\cite{BCG24}.
		
	\subsection{Numerical results}
	
	We now present some numerical tests for fluid-structure interaction problems in two dimensions. First, we focus on a steady-state problem to analyze how the choice of assembly technique for the coupling terms affects the convergence of the method, then we solve two time dependent problems.
	
	\subsubsection{Convergence analysis for stationary problem}
	
	We consider a stationary version of Problem~\ref{pro:discrete} and we study how the convergence of the method is affected by the assembly technique for the coupling term. Here we only discuss the example of the immersed square presented in~\cite{BCG24}. The interested reader can refer to~\cite{boffi2022interface,BCG24} for a wider set of numerical tests.
	
\begin{figure}
		
		\includegraphics[width=0.24\linewidth]{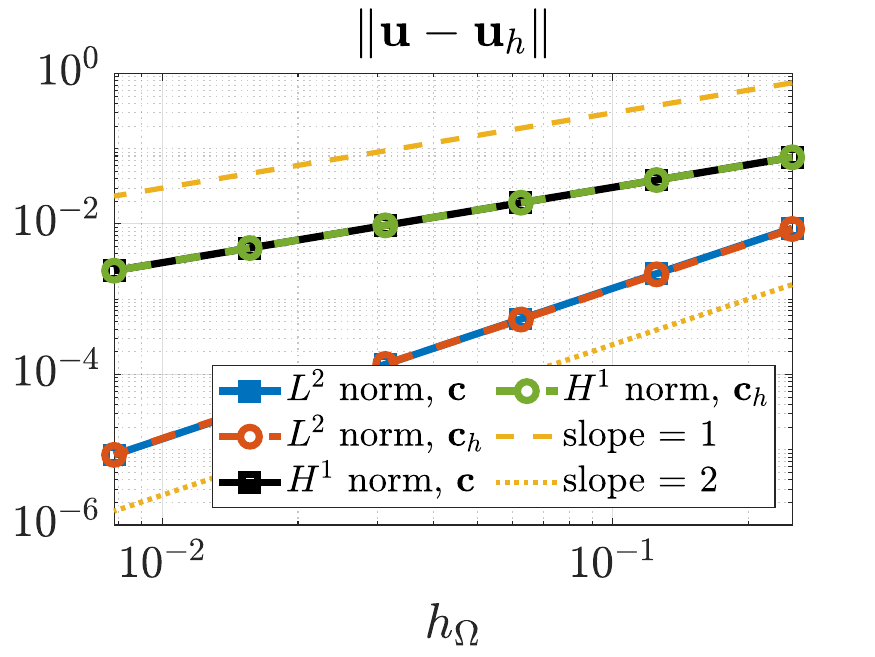}
		\includegraphics[width=0.24\linewidth]{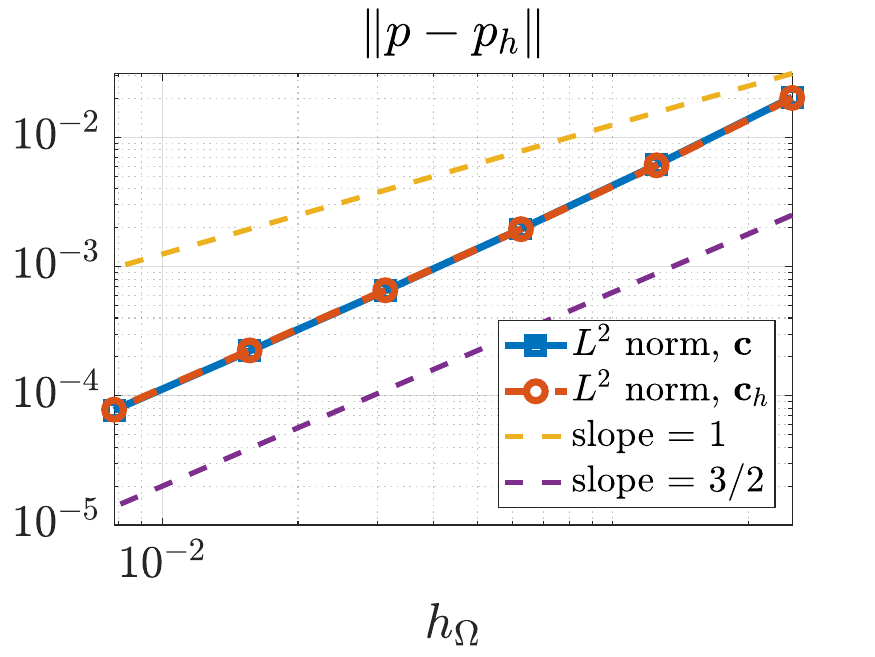}
		\includegraphics[width=0.24\linewidth]{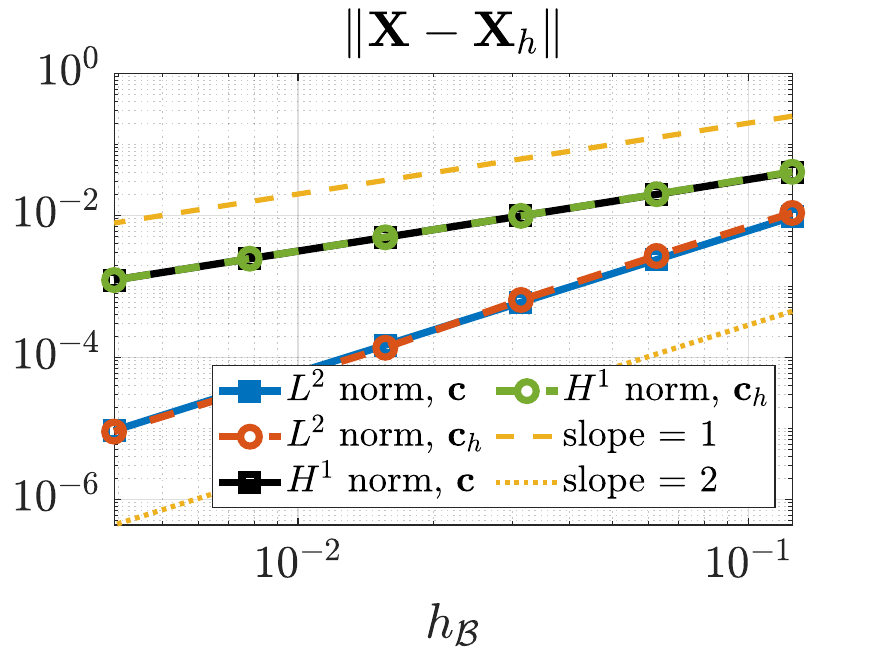}
		\includegraphics[width=0.24\linewidth]{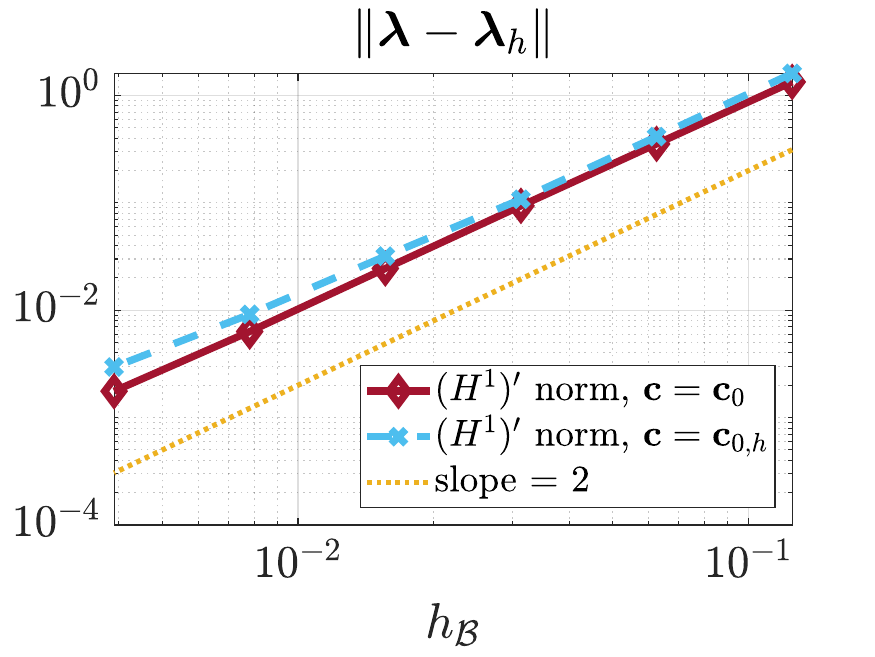}
		
		\caption{Convergence plots for the immersed square problem. The coupling bilinear form is the scalar product in $\LdBd$ computed by exact ($\c$, solid lines) and inexact ($\c_h$, dashed lines) quadrature formulas. Meshes are refined in such a way that $h_\B/h_\Omega=1$.}
		\label{fig:comp1b}
	\end{figure}

	The fluid domain is $\Omega=[-2,2]^2$, whereas the reference and actual configurations of the immersed square are $\B=[0,1]^2$ and $\Os=[-0.62,1.38]^2$, respectively. We compute the right hand side of our problem so that the following exact solution is obtained:
	\begin{equation}\label{eq:solution}
		\begin{aligned}
		&\u(x,y) = \curl\left((4-x^2)^2(4-y^2)^2\right),&&\quad
		p(x,y) = 150\,\sin(x),\\
		&\X(s_1,s_2) = \curl\left((4-s_1^2)^2(4-s_2^2)^2\right),
		&&\quad\llambda(x,y) = [\exp(x),\exp(y)].
		\end{aligned}
	\end{equation}
	We solve the problem on a sequence of uniform triangulations by choosing the $(\Pcal_1iso\Pcal_2,\Pcal_1)$ element for the Stokes equation (also known as Bercovier--Pironneau element~\cite{bercovier}), while the deformation and the Lagrange multiplier are discretized by continuous piecewise linear elements $\Pcal_1$. We compare the behavior of the two choices of coupling bilinear form with respect to the assembly technique.
	Figure~\ref{fig:comp1b} reports the convergence history of the error in the case of the $\LdBd$ coupling term. It is evident that we obtain optimal results for both exact and inexact assembly procedures (solid and dashed lines, respectively).
	
\begin{figure}
		
		\includegraphics[width=0.24\linewidth]{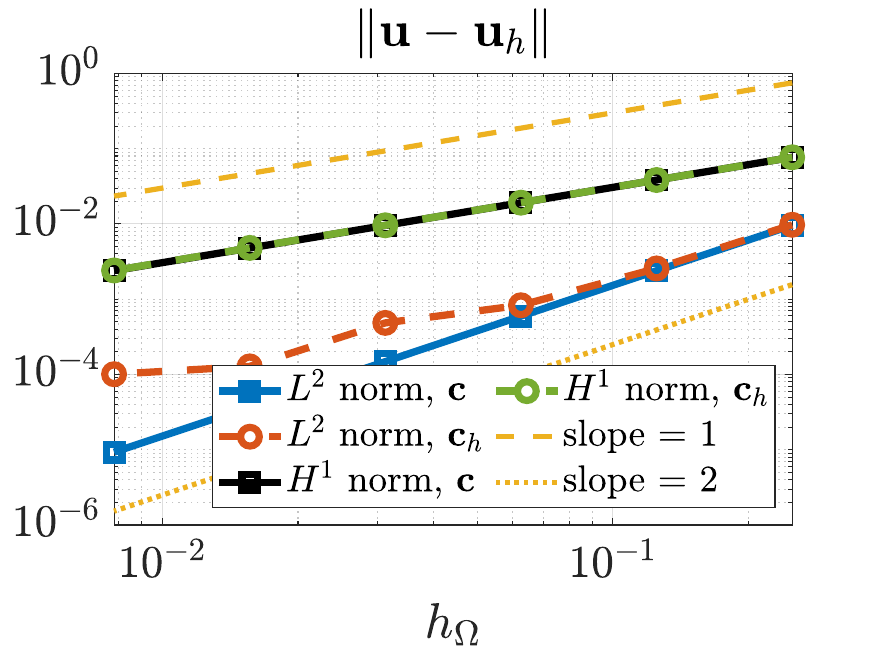}
		\includegraphics[width=0.24\linewidth]{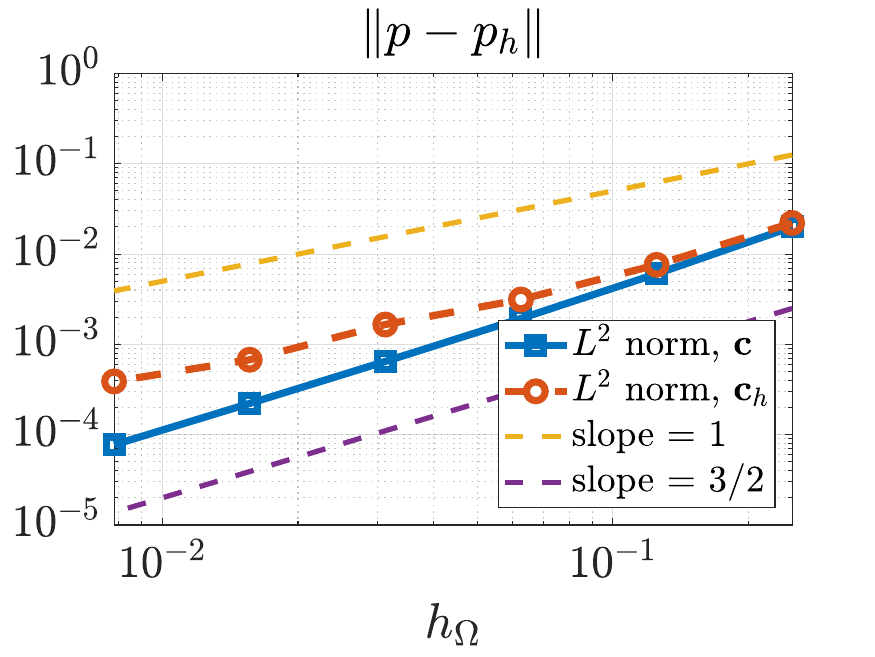}
		\includegraphics[width=0.24\linewidth]{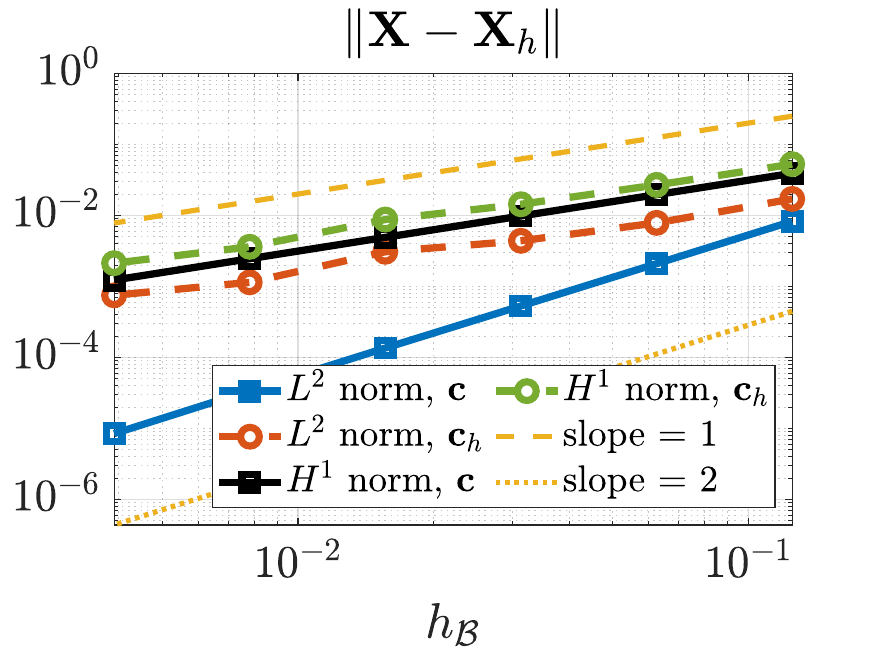}
		\includegraphics[width=0.24\linewidth]{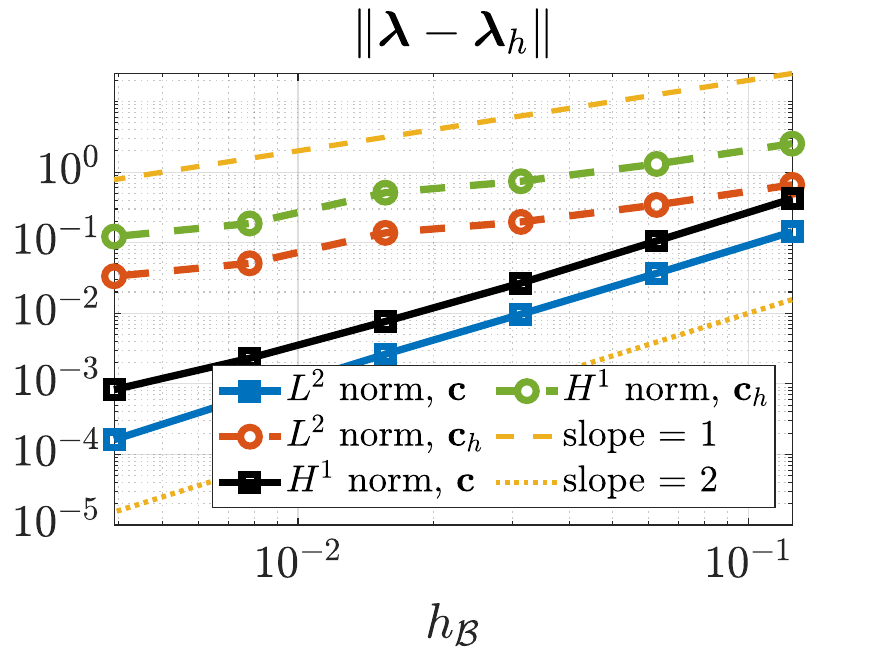}
		
		\caption{Convergence plots for the immersed square problem. The coupling bilinear form is the scalar product in $\Hub$ computed by exact ($\c$, solid lines) and inexact ($\c_h$, dashed lines) quadrature formulas. Meshes are refined in such a way that $h_\B/h_\Omega=1$.}
		\label{fig:comp1a}
		
	\end{figure}

\begin{figure}
		\includegraphics[width=0.24\linewidth]{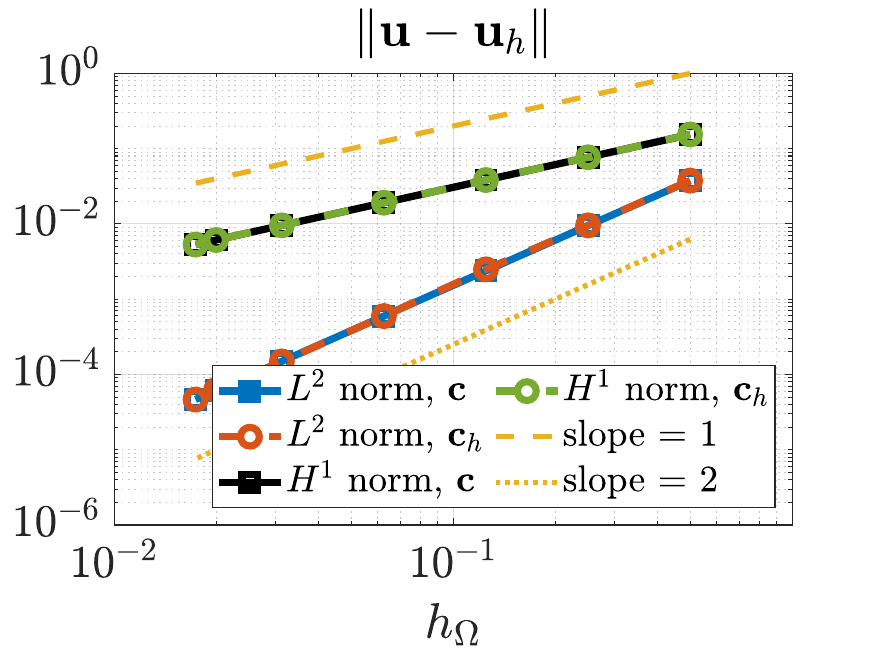}
		\includegraphics[width=0.24\linewidth]{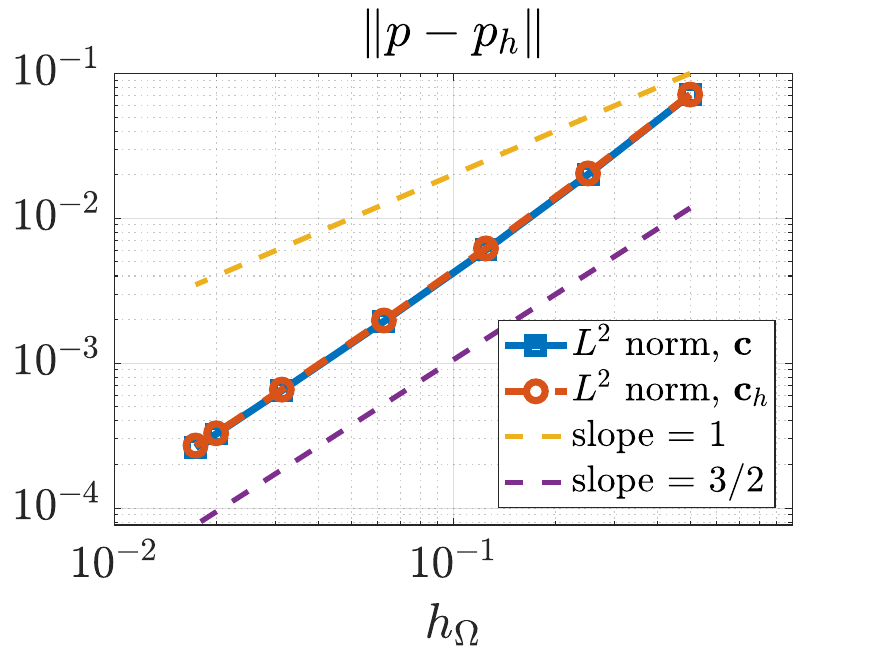}
		\includegraphics[width=0.24\linewidth]{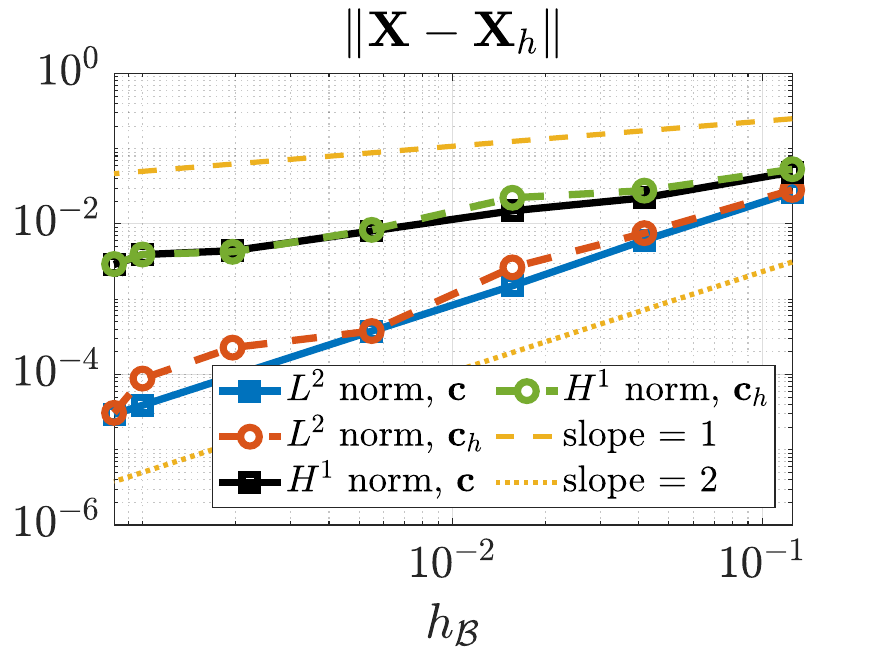}
		\includegraphics[width=0.24\linewidth]{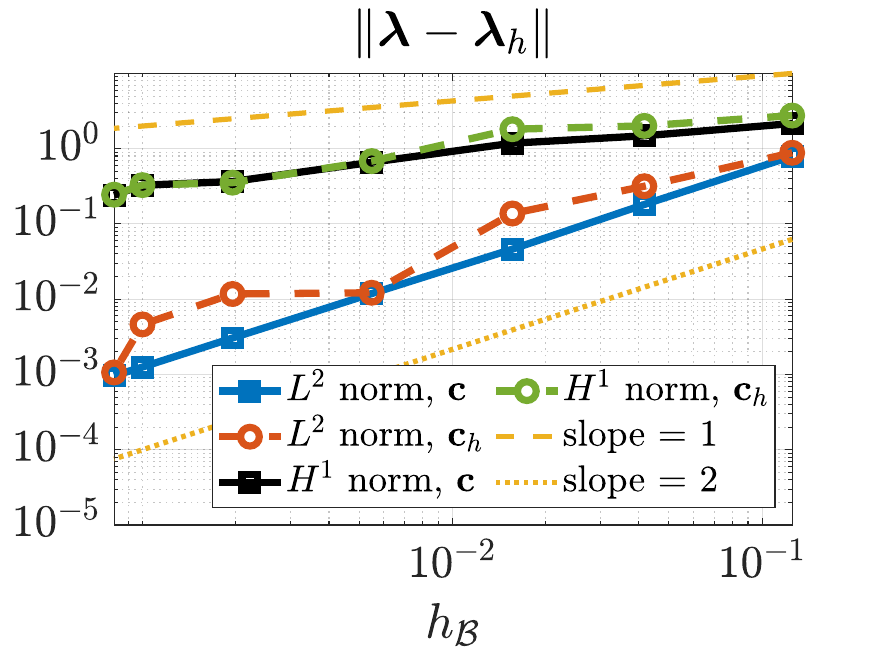}
		
		\caption{Convergence plots for the immersed square problem. The coupling bilinear form is the scalar product in $\Hub$ computed by exact ($\c$, solid lines) and inexact ($\c_h$, dashed lines) quadrature formulas. Meshes are refined in such a way that $h_\B=(h_\Omega/2)^{3/2}$.}
		\label{fig:comp2}
	\end{figure}

	 Figures~\ref{fig:comp1a} and~\ref{fig:comp2} show the convergence history of the error in the case of the $\Hub$ coupling term. In the case of a refinement of fluid and solid meshes with a fixed ratio between the mesh sizes the error deteriorates when the coupling term is computed inexactly (dashed lines). This phenomenon can be circumvented by choosing the mesh sizes so that $h_\B/h_\Omega$ tends to zero. The results reported in Figure~\ref{fig:comp2} are obtained by choosing $h_\B=(h_\Omega/2)^{3/2}$ and provide convergence even if the coupling term is computed in approximate way. Notice that this choice yields a reduced convergence rate. These results are in agreement with the theory presented in~\cite{BCG24}.

	\subsubsection{The stretched annulus}
	
	We consider an annular structure immersed in a steady state Newtonian fluid governed by the Stokes equation. At initial time, the annulus is stretched and, during the simulation, internal elastic forces drive it to its resting configuration, generating the motion of the fluid. 
	This two-dimensional example, already discussed in~\cite{2015,boffi2022parallel}, can be interpreted as the section of a thick cylinder immersed in a rigid box with square section. 
	
	The solid material is modeled by the linear constitutive law $\PP(\F)=\kappa\F$, with $\kappa=10$. At rest, the annulus occupies the region $\{\x\in\R^2:\,0.3\le|\x|\le0.5\}$, being at the center of the box~$[-1,1]^2$ filled by the fluid. Since the geometry of the system is symmetric, we can reduce our study to a quarter of the domain, thus we set
	\begin{equation}
		\Omega=[0,1]^2
		\qquad\text{and}\qquad
		\B=\{\s\in\R^2:\,s_1,s_2\ge0,\,0.3\le|\s|\le0.5\},
	\end{equation}
	with $\s=(s_1,s_2)$.
	The stretching of the annulus is given by the following initial condition for the deformation $\X$
	\begin{equation}
		\X(\s,0) = \left(\frac{s_1}{1.4},1.4\,s_2\right),
	\end{equation}
	while the initial velocity $\u(\x,0)$ is set to zero. Moreover, we impose no slip boundary conditions for $\u$ on the upper and right edge of $\Omega$, while both fluid and structure can move along the bottom and left edges in tangential direction.
	
	We also assume that fluid and solid have same viscosity: this is a reasonable mathematical assumption when dealing with biological systems, see~\cite{peskin1972flow} for instance. More precisely, $\nu=0.1$. We mentioned in the introduction that the added mass effect arises when fluid and solid have equal or similar densities. In order to confirm the robustness of our method, we set $\rho_f=\rho_s=1$.
	
	The equations are discretized by the $(\mathcal{Q}_2,\Pcal_1,\mathcal{Q}_1,\mathcal{Q}_1)$ finite element. The linear system is solved using GMRES accelerated by the preconditioners introduced in Section~\ref{sec:solvers}, whose action is performed through exact inversion. The performance of the parallel solver has been analyzed, in terms of optimality and scalability, in~\cite{boffi2022parallel}. The coupling term is the scalar product in $\LdB$ assembled exactly.
	
\begin{figure}
		\begin{center}
			\includegraphics[trim = 30 10 20 0,clip, width=3.5cm]{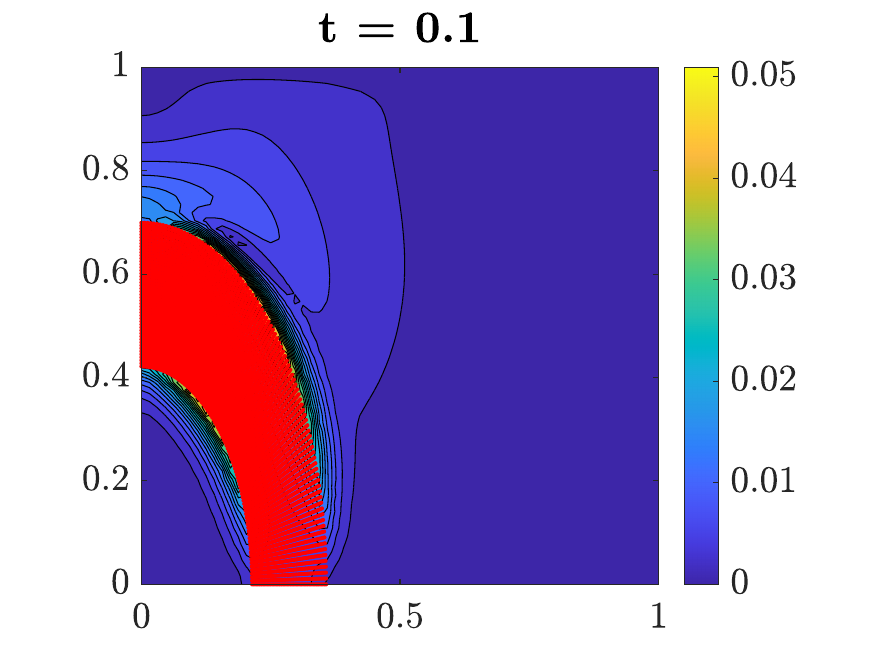}
			\includegraphics[trim = 30 10 20 0,clip, width=3.5cm]{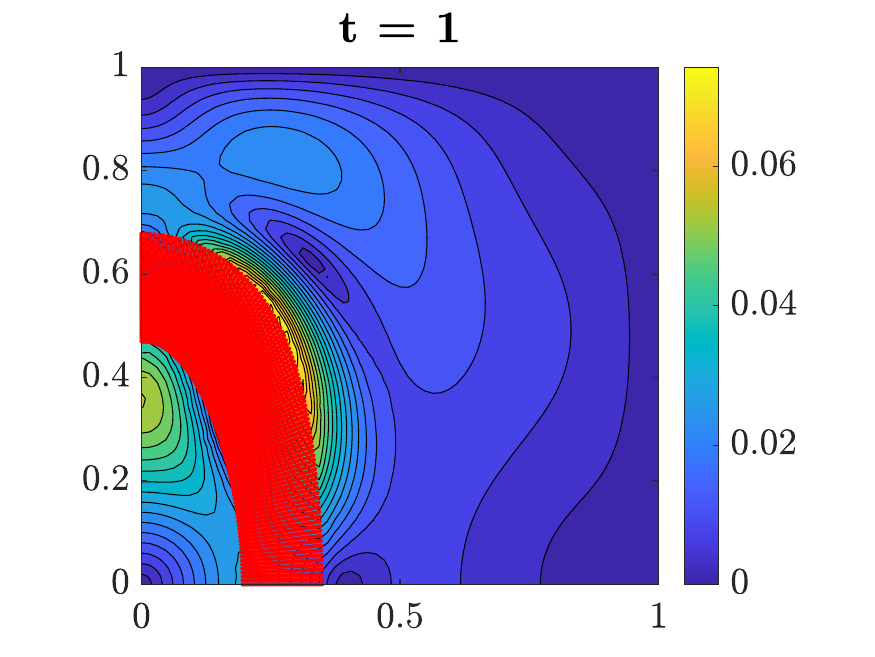}
			\includegraphics[trim = 30 10 20 0,clip, width=3.5cm]{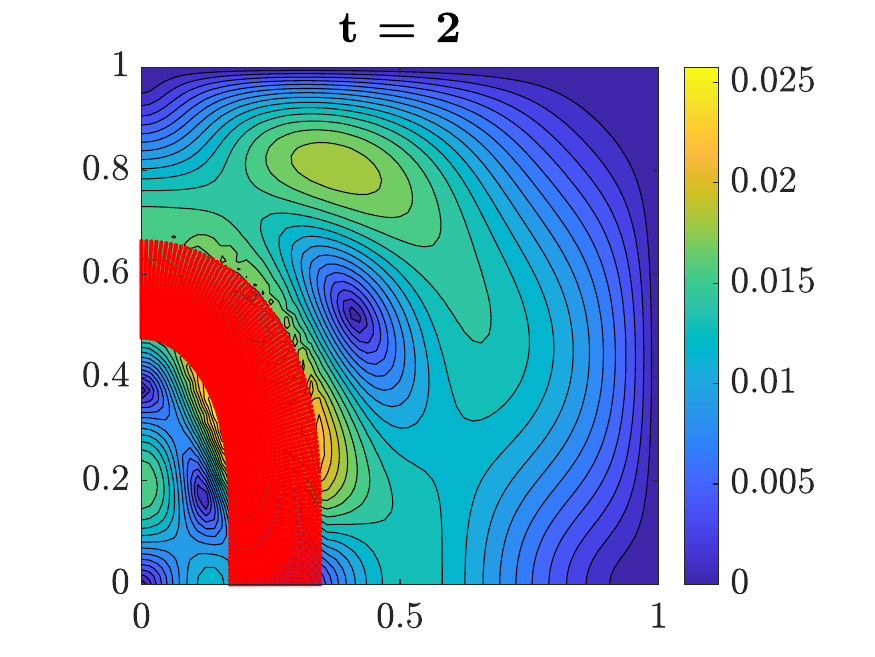}
			\includegraphics[trim = 30 10 20 0,clip, width=3.5cm]{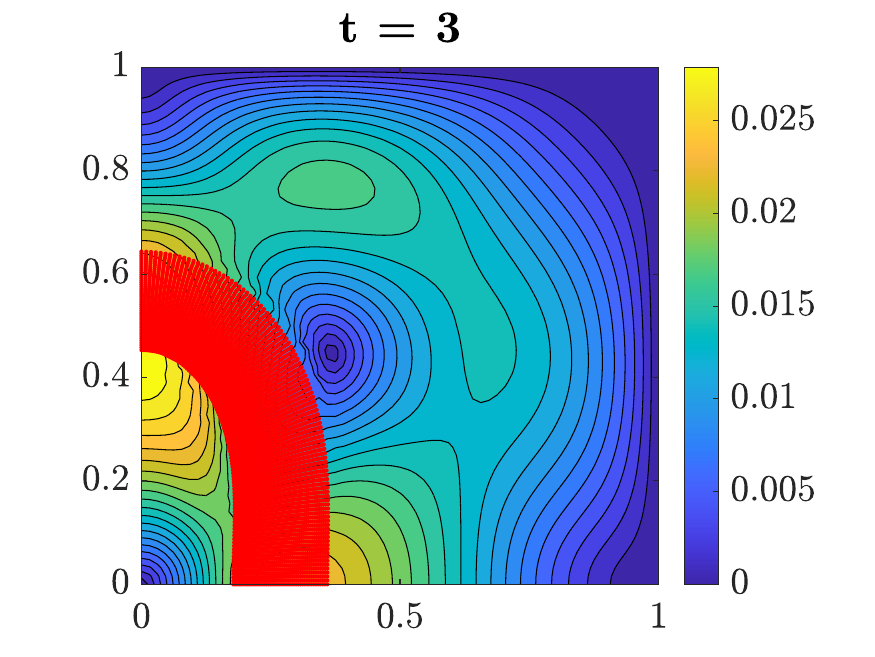}\\
			\medskip
			\includegraphics[trim = 30 10 20 0,clip, width=3.5cm]{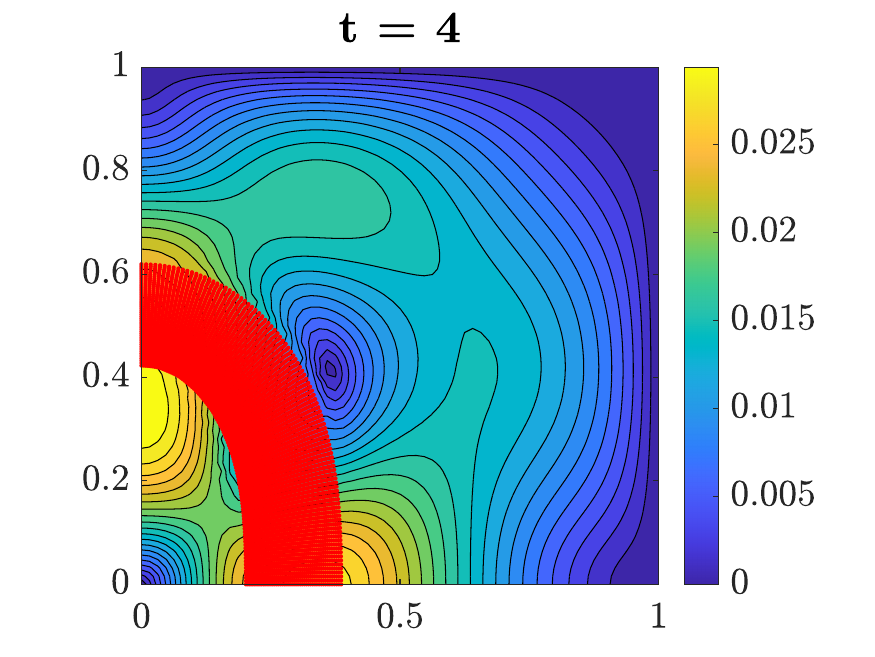}
			\includegraphics[trim = 30 10 20 0,clip, width=3.5cm]{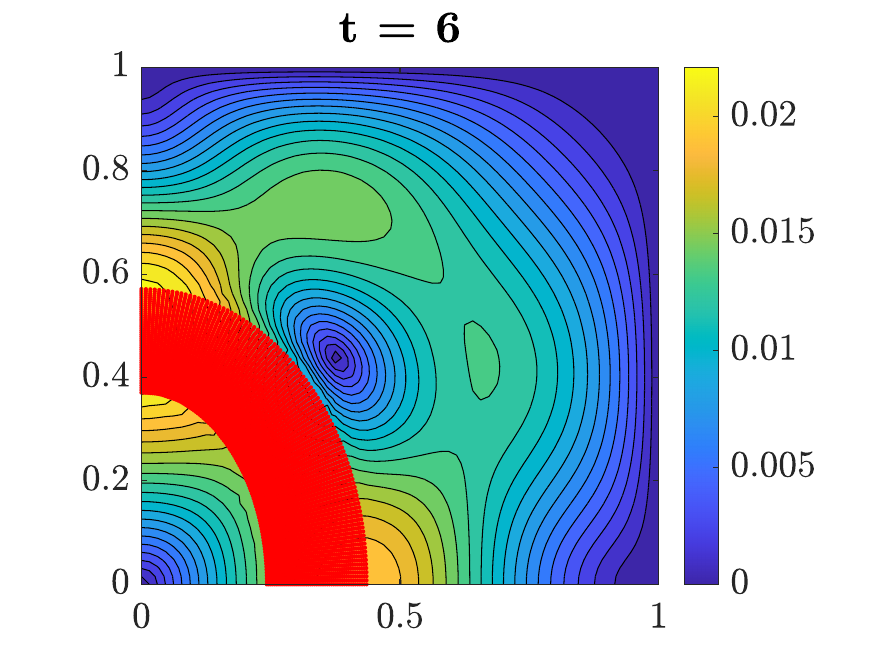}
			\includegraphics[trim = 30 10 20 0,clip, width=3.5cm]{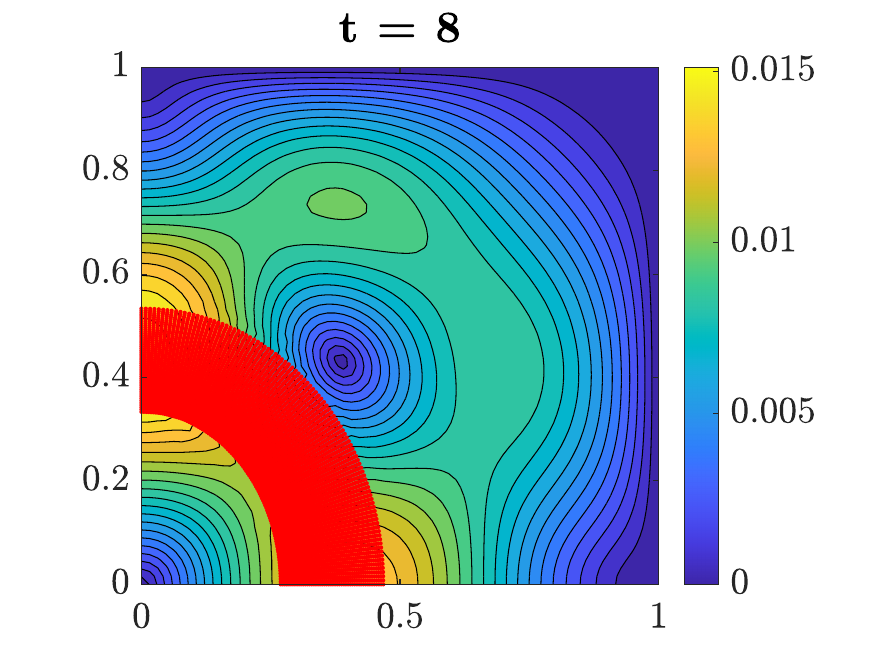}
			\includegraphics[trim = 30 10 20 0,clip, width=3.5cm]{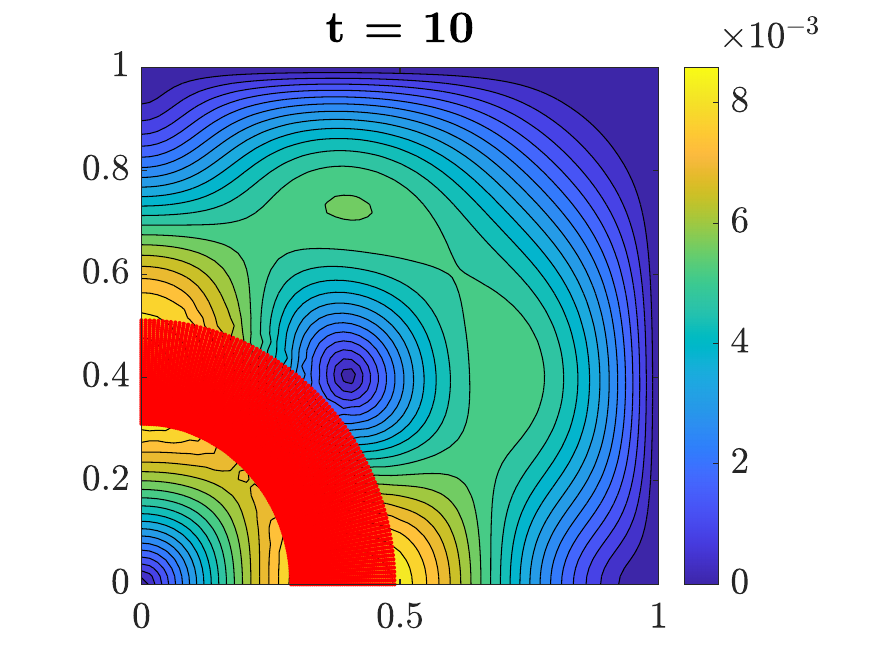}\\
			\medskip
			\includegraphics[trim = 30 10 20 0,clip, width=3.5cm]{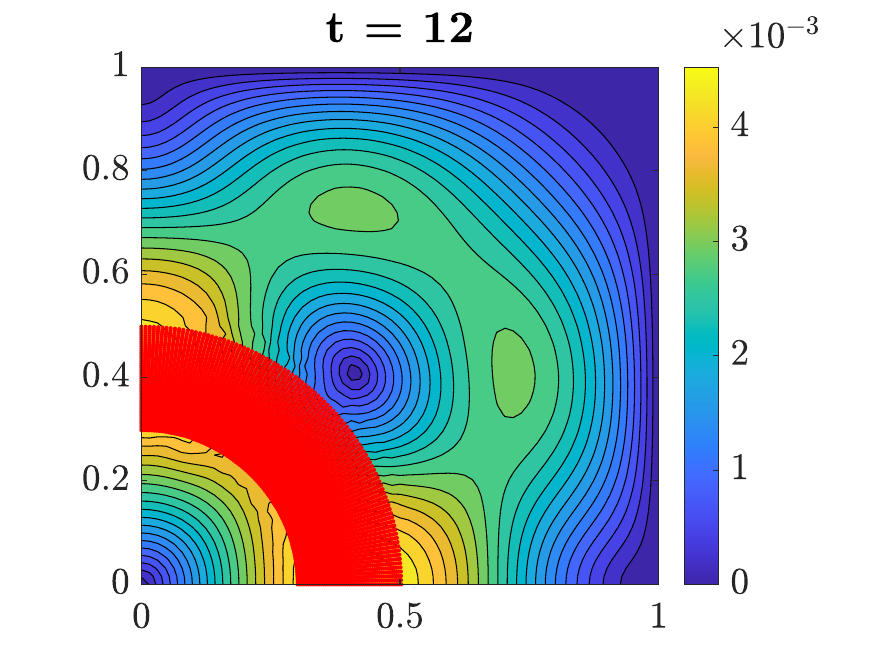}
			\includegraphics[trim = 30 10 20 0,clip, width=3.5cm]{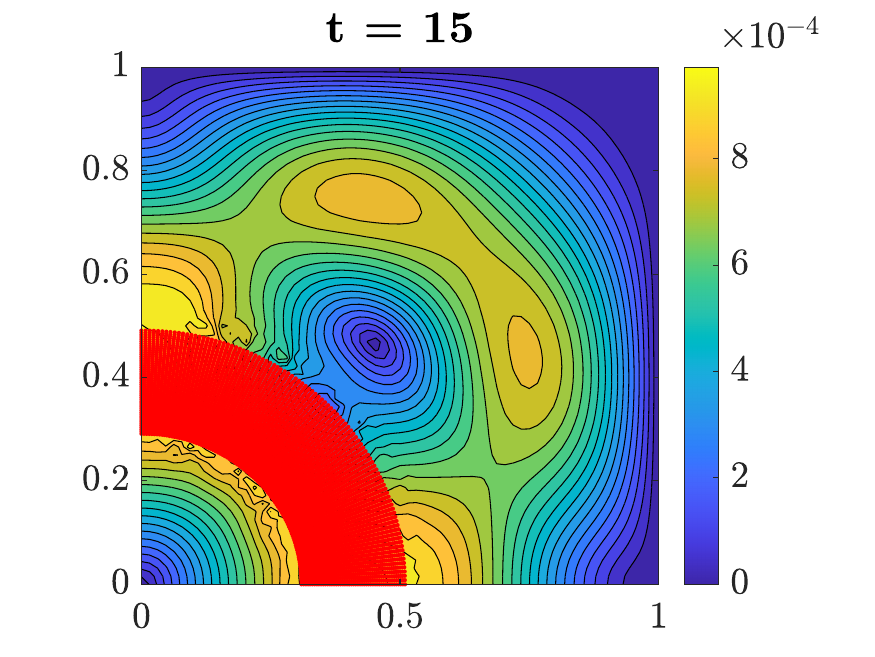}
			\includegraphics[trim = 30 10 20 0,clip, width=3.5cm]{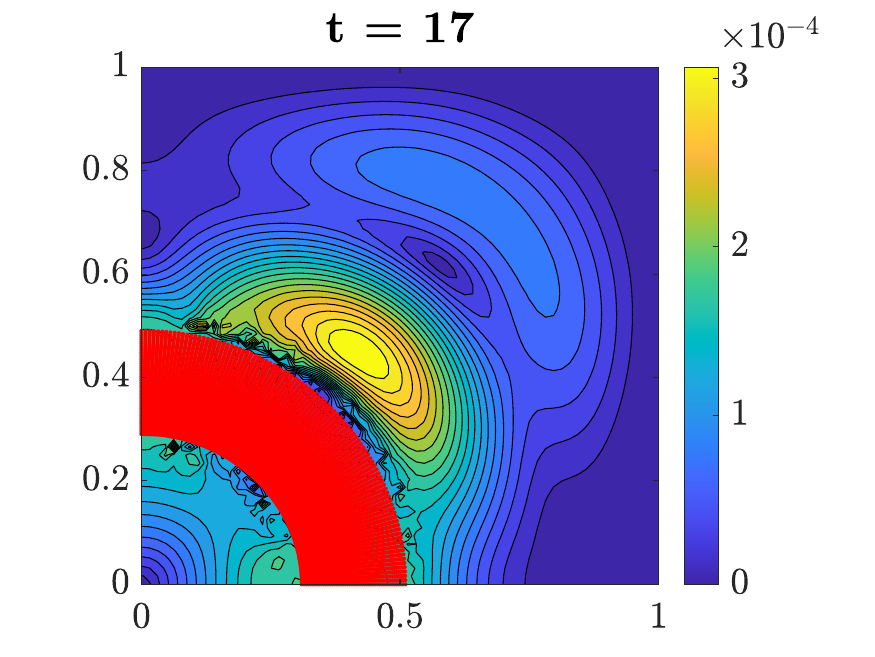}
			\includegraphics[trim = 30 10 20 0,clip, width=3.5cm]{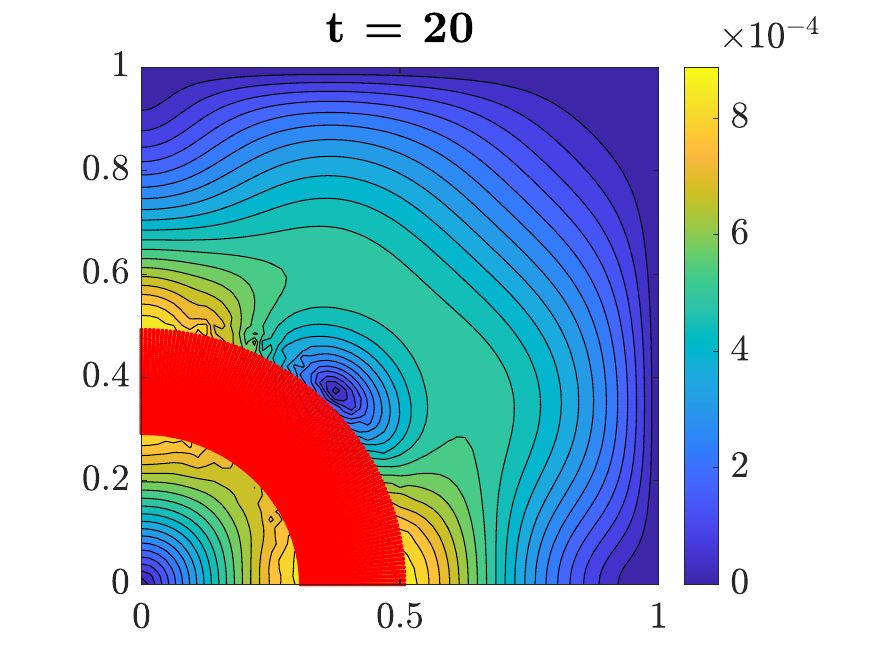}\\
			\caption{Snapshots of the annulus evolution (red) with velocity streamlines (colorbar).}
			\label{fig:annulus}
		\end{center}
	\end{figure}

	Some snapshots of the simulation are reported in Figure~\ref{fig:annulus}: we solved the problem on a mesh with 30,534 total degrees of freedom until time $T=20$ with step $\dt=0.01$. The solver is stable and the results are consistent with the physics of the problem.
		
	
	
	\subsubsection{The bar}
	
	This test has been presented in~\cite{boffi2022parallel}. We consider the Stokes equation in $\Omega=[0.1]^2$. The immersed structure is an elastic bar anchored at the left edge of $\Omega$. The solid material is modeled by the nonlinear constitutive law given by the exponential strain energy function of an isotropic hyperelastic material, that is
	$$
	W(\F) = \frac{\gamma}{2\eta}\exp\left(\eta\,\mathrm{tr}(\F^\top\F)-2\right),
	$$
	where $\mathrm{tr}(\F^\top\F)$ denotes the trace of $\F^\top\F$, $\gamma=1.333$ and $\eta=9.242$. The reference domain for the bar is $\B=[0,0.4]\times[0.45,0.55]$ and corresponds to the initial configuration, thus $\X(\s,0)=\s$. Moreover, $\u(\X,0)=0$.
	
	During the time interval $[0,1]$ the structure is pulled down by a force applied at the middle point of its right edge. Then the structure is released so that internal forces bring it back to rest.
	
	\begin{figure}
		\begin{center}
			\includegraphics[trim = 30 10 20 0,clip, width=3.5cm]{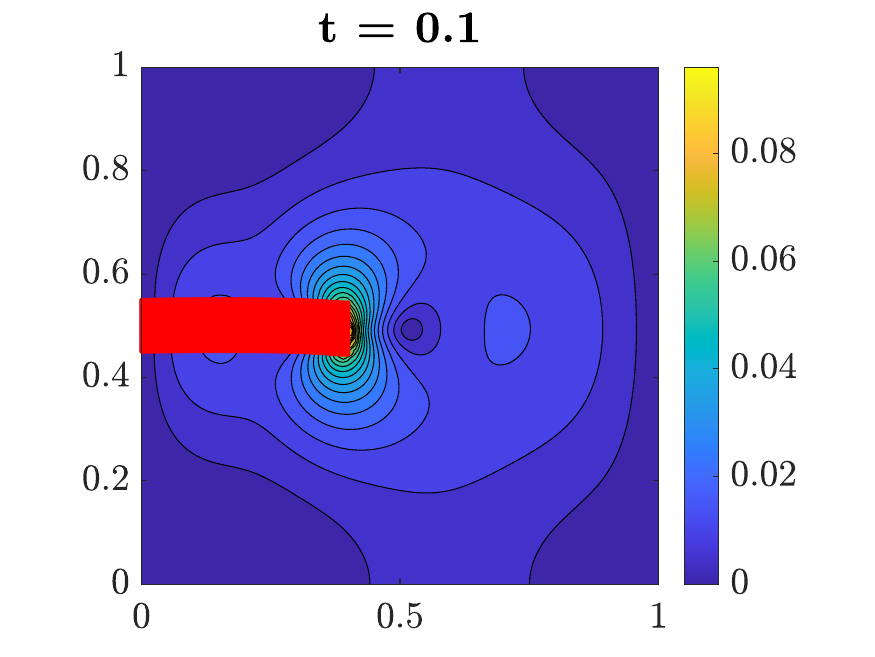}
			\includegraphics[trim = 30 10 20 0,clip, width=3.5cm]{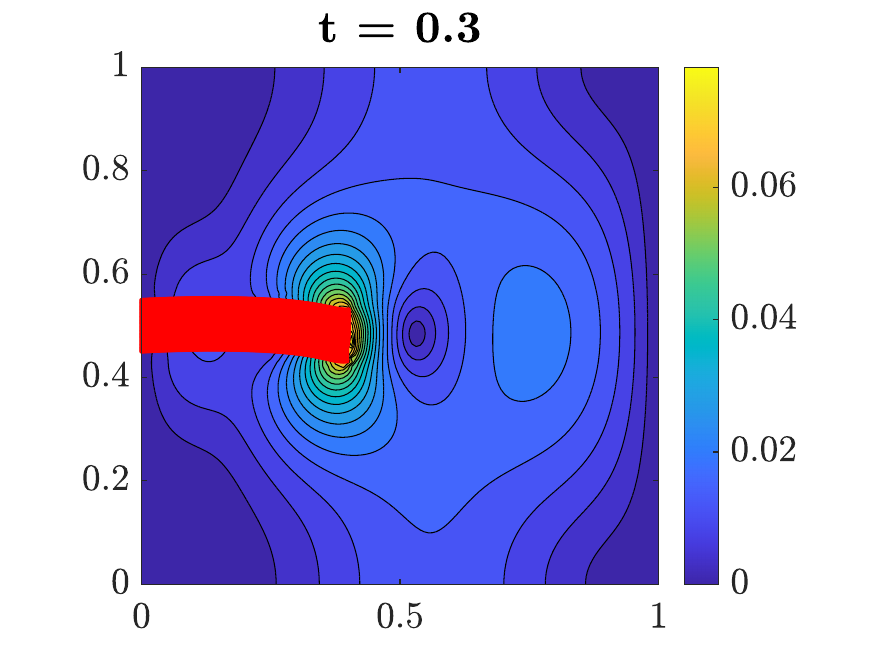}
			\includegraphics[trim = 30 10 20 0,clip, width=3.5cm]{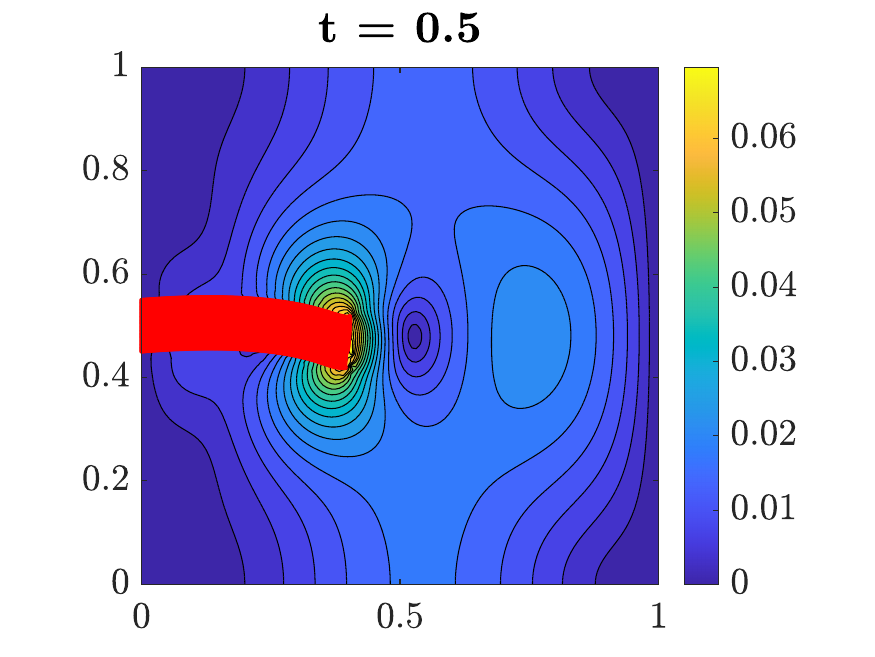}
			\includegraphics[trim = 30 10 20 0,clip, width=3.5cm]{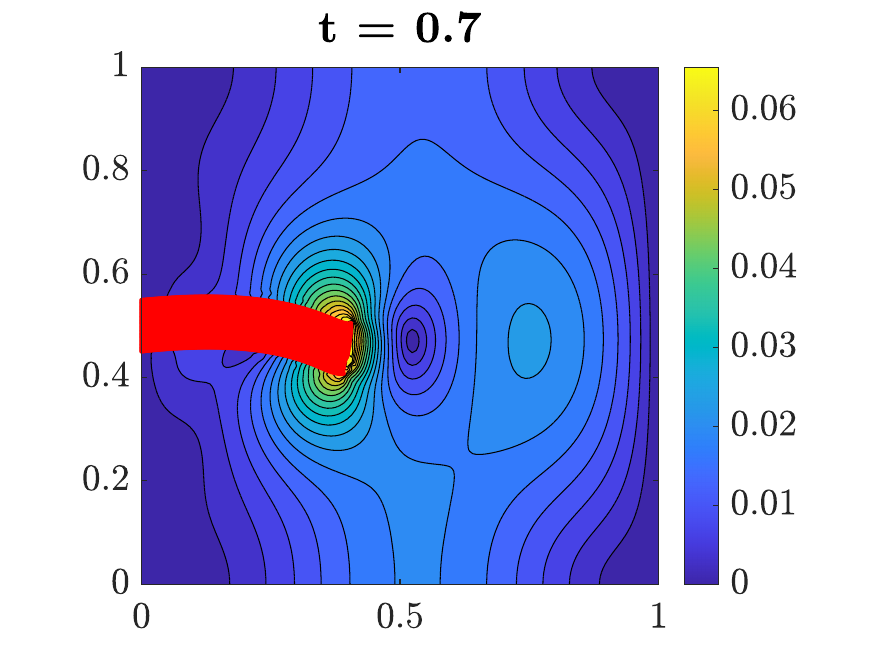}\\
			\medskip
			\includegraphics[trim = 30 10 20 0,clip, width=3.5cm]{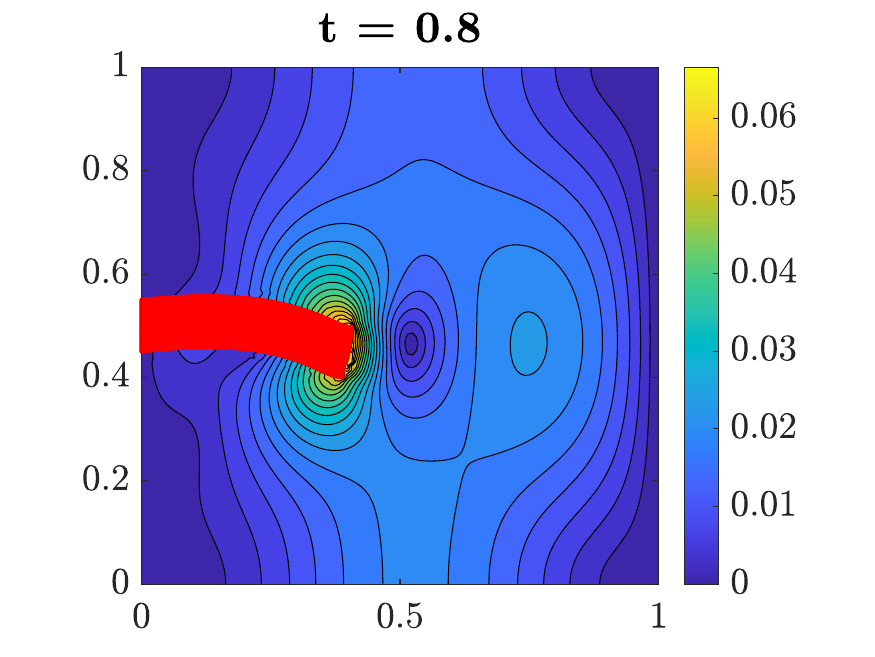}
			\includegraphics[trim = 30 10 20 0,clip, width=3.5cm]{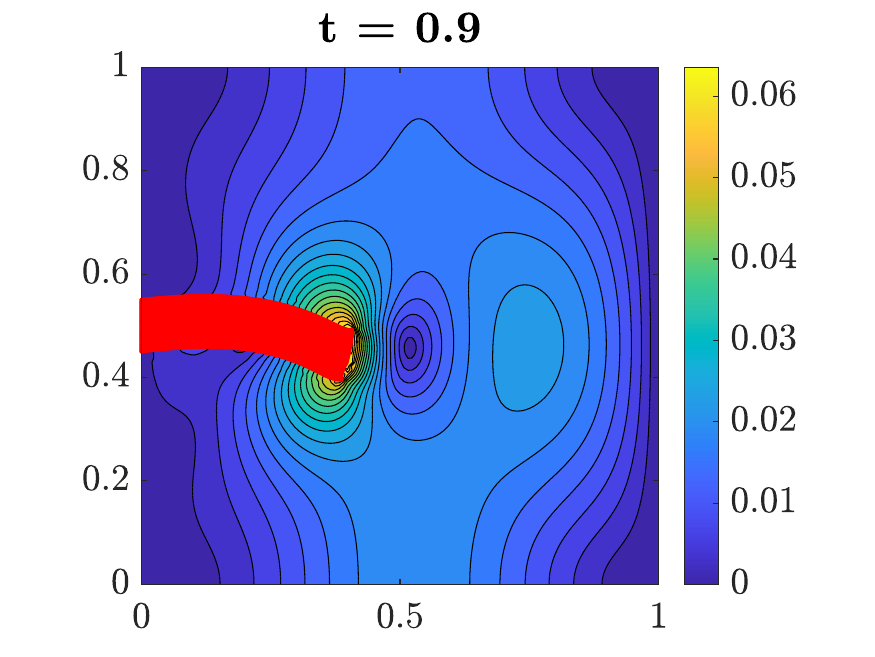}
			\includegraphics[trim = 30 10 20 0,clip, width=3.5cm]{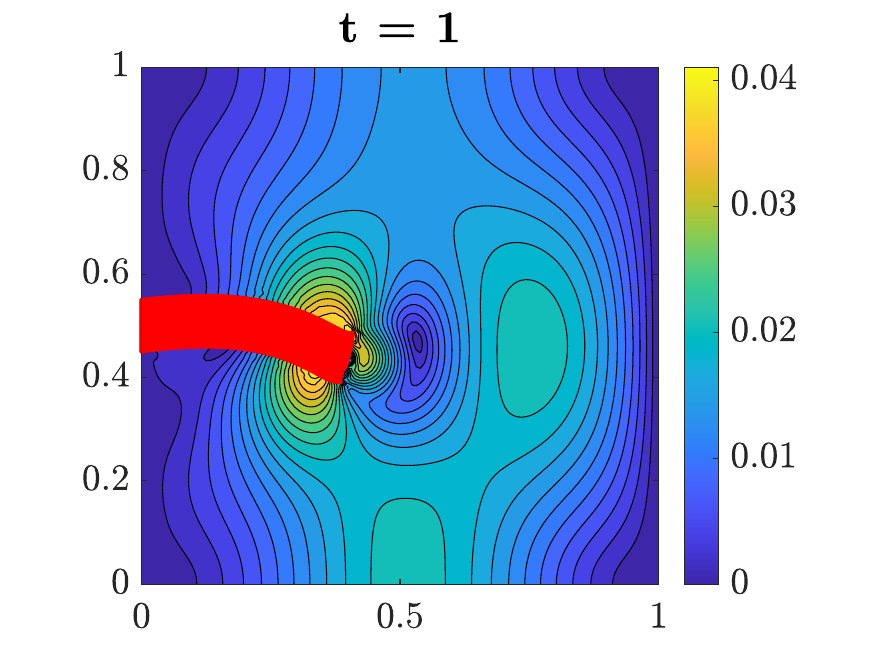}
			\includegraphics[trim = 30 10 20 0,clip, width=3.5cm]{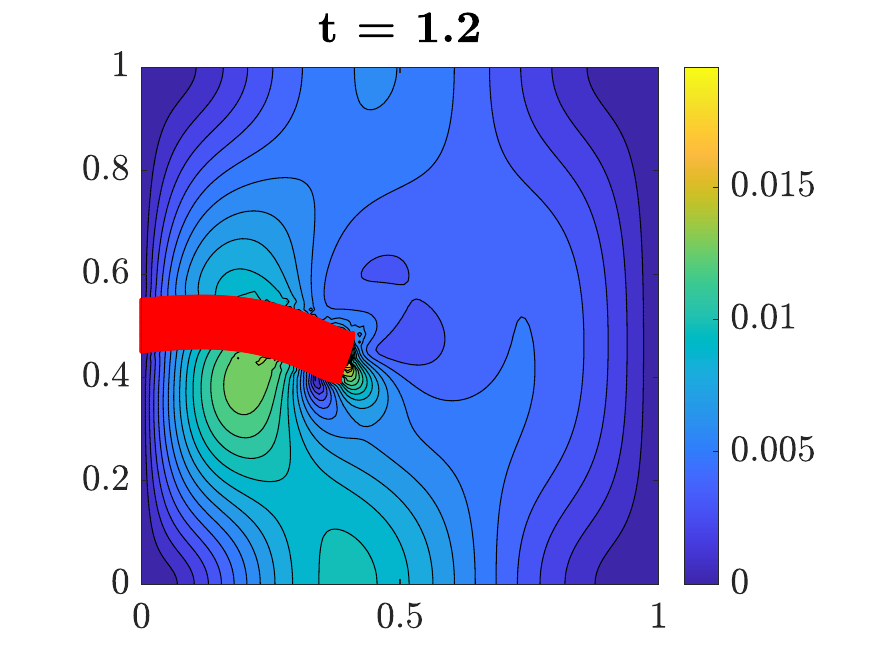}\\
			\medskip
			\includegraphics[trim = 30 10 20 0,clip, width=3.5cm]{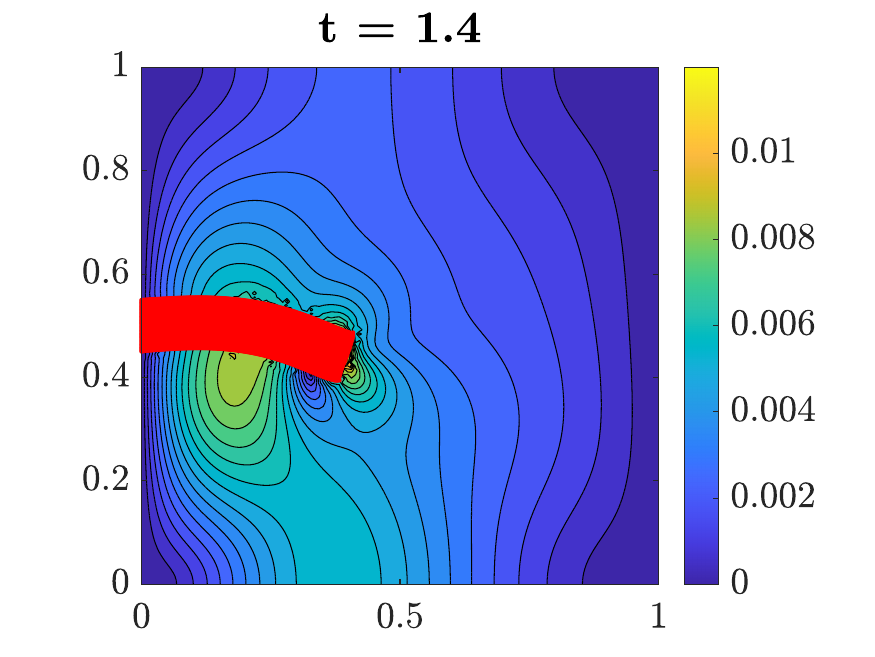}
			\includegraphics[trim = 30 10 20 0,clip, width=3.5cm]{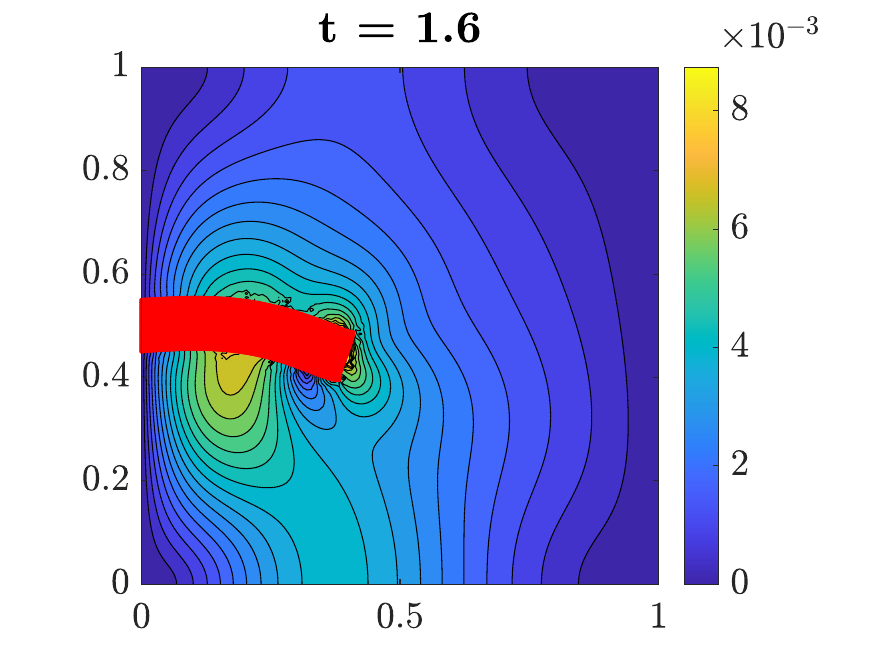}
			\includegraphics[trim = 30 10 20 0,clip, width=3.5cm]{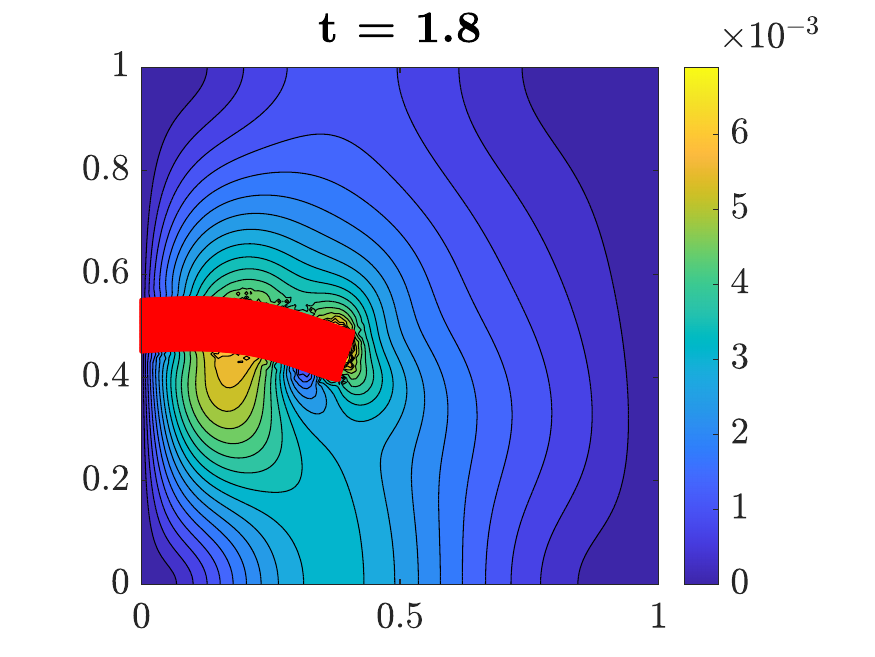}
			\includegraphics[trim = 30 10 20 0,clip, width=3.5cm]{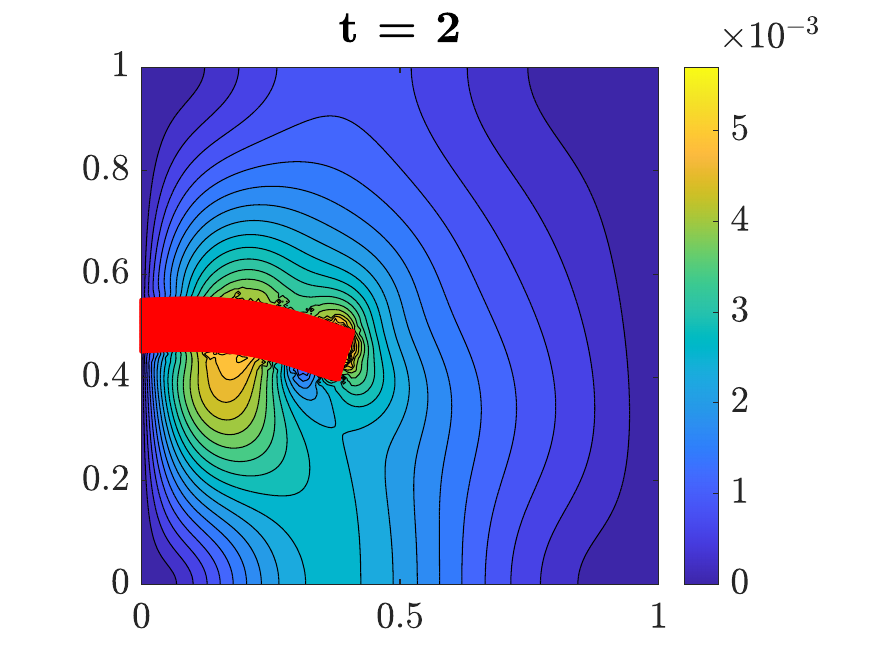}\\
			\caption{Snapshots of the bar evolution (red) with velocity streamlines (colorbar).}
			\label{fig:bar_evolution}
		\end{center}
	\end{figure}
	
	We consider again the  $(\mathcal{Q}_2,\Pcal_1,\mathcal{Q}_1,\mathcal{Q}_1)$ finite element. In order to tackle the nonlinearity of the solid equation, we employ the Newton method. The linear system is then solved by GMRES combined with $\mathsf{P}_{\mathsf{diag}}$ and $\mathsf{P}_{\mathsf{tri}}$, which are inverted exactly. The coupling between fluid and solid is enforced by the scalar product in $\LdBd$ assembled with the exact procedure.
	
	We simulate the system in the time interval $[0,2]$ with $\dt=0.01$ on a sequence of six meshes. The total number of degrees of freedom is 83,398. Some snapshots are depicted in Figure~\ref{fig:bar_evolution}. The optimality and scalability of the parallel solver have been discussed in~\cite{boffi2022parallel}. Also in this case, the solver is stable and the results are consistent with the physics of the problem.

	\section{Conclusions}
	
	We discussed the latest results on an unfitted finite element method for interface and fluid-structure interaction problems based on a Lagrange multiplier in the spirit of the fictitious domain approach. 
	
	Due to jumps in the coefficient or the immersion of the solid body in the fluid, the physical domain is partitioned into two regions. The main idea behind our approach is to extend one domain to the other. Then, the equations are solved on both the background and the immersed domains, while a coupling term describes their interaction in the fictitiously extended region. The computational domains are fixed and discretized by two independent meshes.
	
	After recalling the main stability properties of our approach, we discussed how to deal with the coupling term. Indeed, the computation of its discrete version involves integration of functions defined over two distinct non-matching grids. The coupling term can be computed exactly by considering a composite quadrature rule on the intersection between the two meshes, or in a (cheaper) approximate way taking into account the presence of an additional error.
	
	We then mentioned the main difficulties arising when designing an efficient solver for our formulation. We presented two suitable preconditioners able to accelerate the convergence of common iterative solvers and we compared their performance in terms of solution time.
	
	Finally, we addressed the issue of mesh adaptivity by introducing optimal a posteriori error estimators, satisfying both the reliability and efficiency properties.
	
	The description of computational aspects and theoretical results is completed by some numerical tests on both elliptic interface problems and fluid-structure interactions.
	
	\section*{Acknowledgments}
	\noindent
	D. Boffi, F. Credali and L. Gastaldi are members of GNCS/INdAM Research group. The research of L. Gastaldi is partially supported by PRIN/MUR (grant No.20227K44ME) and IMATI/CNR.

\bibliographystyle{unsrt}
\bibliography{biblio}

\end{document}